\documentclass[tbtags,reqno]{amsart}
\usepackage{amsthm,amsopn,amsfonts,amssymb,epsfig} 

\catcode `!=11

\newdimen\squaresize 
\newdimen\thickness 
\newdimen\Thickness
\newdimen\ll! \newdimen \uu! \newdimen\dd! \newdimen \rr! 
\newdimen \temp!

\def\sq!#1#2#3#4#5{%
\ll!=#1 \uu!=#2 \dd!=#3 \rr!=#4
\setbox0=\hbox{%
 \temp!=\squaresize\advance\temp! by .5\uu!
 \rlap{\kern -.5\ll! 
 \vbox{\hrule height \temp! width#1 depth .5\dd!}}%
 \temp!=\squaresize\advance\temp! by -.5\uu!  
 \rlap{\raise\temp! 
 \vbox{\hrule height #2 width \squaresize}}%
 \rlap{\raise -.5\dd!
 \vbox{\hrule height #3 width \squaresize}}%
 \temp!=\squaresize\advance\temp! by .5\uu!
 \rlap{\kern \squaresize \kern-.5\rr! 
 \vbox{\hrule height \temp! width#4 depth .5\dd!}}%
 \rlap{\kern .5\squaresize\raise .5\squaresize
 \vbox to 0pt{\vss\hbox to 0pt{\hss $#5$\hss}\vss}}%
}
 \ht0=0pt \dp0=0pt \box0
}

\def\vsq!#1#2#3#4#5\endvsq!{\vbox to \squaresize{\hrule
width\squaresize height 0pt%
\vss\sq!{#1}{#2}{#3}{#4}{#5}}}

\newdimen \LL! \newdimen \UU! \newdimen \DD! \newdimen \RR!

\def\vvsq!{\futurelet\next\vvvsq!}
\def\vvvsq!{\relax
  \ifx     \next l\LL!=\Thickness \let\continue=\skipnexttoken!
  \else\ifx\next u\UU!=\Thickness \let\continue=\skipnexttoken!
  \else\ifx\next d\DD!=\Thickness \let\continue=\skipnexttoken!
  \else\ifx\next r\RR!=\Thickness \let\continue=\skipnexttoken!
  \else\def\continue{\vsq!\LL!\UU!\DD!\RR!}%
  \fi\fi\fi\fi
  \continue}

\def\skipnexttoken!#1{\vvsq!}

\def\place#1#2#3{\vbox to 0pt{\vss
\rlap{\kern#1\squaresize
  \raise#2\squaresize\hbox{$#3$}}
\vss}}

\def\Young#1{\LL!=\thickness \UU!=\thickness \DD! = \thickness \RR! =
\thickness
\vbox{\smallskip\offinterlineskip
\halign{&\vvsq! ## \endvsq!\cr #1}}}

\def\blank{\omit\hskip\squaresize}
\catcode `!=12

\edef\savecatcodeat{\the\catcode`@}
\catcode`\@=11

\def\tb@ifSpecChars#1#2{#1}
\def\tb@ifNoSpecChars#1#2{#2}

\def\tableau{%
  \bgroup
  \@ifstar{\let\Tif\tb@ifNoSpecChars\tb@tableauB}
          {\let\Tif\tb@ifSpecChars\tb@tableauB}}

\def\tb@tableauB{
  \@ifnextchar[{\tb@tableauC}{\tb@tableauC[]}}

\def\tb@tableauC[#1]{\hbox\bgroup%
    \let\\=\cr
    \def\bl{\global\let\tbcellF\tb@cellNF}%
    \def\tf{\global\let\tbcellF\tb@cellH}
    \dimen2=\ht\strutbox \advance\dimen2 by\dp\strutbox%
    \ifx\baselinestretch\undefined\relax%
    \else%
       \dimen0=100sp \dimen0=\baselinestretch\dimen0%
       \dimen2=100\dimen2 \divide\dimen2 by\dimen0%
    \fi%
    \let\tpos\tb@vcenter
    \tb@initYoung
    \tb@options#1\eoo
    \let\arrow\tb@arrow%
    \dimen0=\Tscale\dimen2%
    \dimen1=\dimen0 \advance\dimen1 by \tb@fframe%
    \lineskip=0pt\baselineskip=0pt
    \def\tb@nothing{}%
    \def\endcellno{$\rss\egroup\bss\egroup}
    \def\endcell{\endcellno\kern-\dimen0}
    \def\begincell{\vbox to\dimen0\bgroup\vss\hbox to\dimen0\bgroup\hss$}%
    \let\overlay\tb@overlay%
    \let\fl\tb@fl%
    \let\lss\hss\let\rss\hss\let\tss\vss\let\bss\vss
    \def\mkcell##1{
        \let\tbcellF\tb@cellD
        \def\tb@cellarg{##1}
        \ifx\tb@cellarg\tb@nothing\let\tb@cellarg\tb@cellE\fi%
        \begincell\tb@cellarg\endcellno
        \tbcellF}
    \let\savecellF\tbcellF
     \Tif{\catcode`,=4\catcode`|=\active}{}\tb@tableauD}%

\let\tb@savetableauD\tableauD
{
    \catcode`|=\active \catcode`*=\active \catcode`~=\active%
    \catcode`@=\active
\gdef\tableauD#1{%
  \Tif{
    \mathcode`|="8000 \mathcode`*="8000%
    \mathcode`~="8000 \mathcode`@="8000%
    \def@{\bullet}%
    \let|\cr
    \let*\tf
    \let~\sk
  }{}%
  \tpos{\tabskip=0pt\halign{&\mkcell{##}\cr#1\crcr}}%
  \global\let\tbcellF\savecellF
  \egroup
  \egroup}
}
\let\tb@tableauD\tableauD
\let\tableauD\tb@savetableauD
\let\tb@savetableauD\undefined

\def\tb@options#1{\ifx#1\eoo\relax\else\tb@option#1\expandafter\tb@options\fi}

\def\tb@option#1{%
  \if#1t\let\tpos\tb@vtop\fi
  \if#1c\let\tpos\tb@vcenter\fi
  \if#1b\let\tpos\vbox\fi
  \if#1F\tb@initFerrers\fi
  \if#1Y\tb@initYoung\fi
  \if#1s\tb@initSmall\fi
  \if#1m\tb@initMedium\fi
  \if#1l\tb@initLarge\fi
  \if#1p\tb@initPartition\fi
  \if#1a\tb@initArrow\fi
}

\def\tb@vcenter#1{\ifmmode\vcenter{#1}\else$\vcenter{#1}$\fi}

\def\tb@vtop#1{\hbox{\raise\ht\strutbox\hbox{\lower\dimen0\vtop{#1}}}}

\def\tb@initPartition{\def\Tscale{.3}}
\def\tb@initSmall{\def\Tscale{1}}
\def\tb@initMedium{\def\Tscale{2}}
\def\tb@initLarge{\def\Tscale{3}}

\def\tb@initArrow{\dimen2=1.25em}

\def\tb@initYoung{%
  \def\tb@cellE{}
  \let\tb@cellD\tb@cellN
  \def\sk{\global\let\tbcellF\tb@cellNF}}
\def\tb@initFerrers{%
  \def\tb@cellE{\bullet}
  \let\tb@cellD\tb@cellNF
  \def\sk{\bullet}}

\tb@initMedium

\def\tb@sframe#1{%
  \vbox to0pt{
    \vss
    \hbox to0pt{%
      \hss
      \vbox to\dimen1{
        \hrule depth #1 height0pt
        \vss
        \hbox to\dimen1{
          \vrule width #1 height\dimen1
          \hss
          \vrule width #1
          }%
        \vss
        \hrule height #1 depth 0in
        }%
      \kern-\tb@hframe
      }%
    \kern-\tb@hframe}}

\def\tb@hframe{.2pt}\def\tb@fframe{.4pt}\def\tb@bframe{1.2pt}
\def\tb@cellH{\tb@sframe{\tb@bframe}}       
\def\tb@cellNF{}                            
\def\tb@cellN{\tb@sframe{\tb@fframe}}       
\let\tbcellF\tb@cellN                       

\def\tb@rpad{1pt}
\def\tb@lpad{1pt}
\def\tb@tpad{1.8pt}
\def\tb@bpad{1.8pt}

\def\tb@overlay{\endcell\@ifnextchar[{\tb@overlaya}{\begincell}}
\def\tb@overlaya[#1]{\vbox to\dimen0\bgroup%
  \tb@overlayoptions#1\eoo%
  \tss\hbox to\dimen0\bgroup\lss}
\def\tb@overlayoptions#1{\ifx#1\eoo\relax\else\tb@overlayoption#1\expandafter\tb@overlayoptions\fi}

\def\tb@overlayoption#1{
  \if#1t\def\tss{\vskip\tb@tpad}\let\bss\vss\fi
  \if#1c\let\tss\vss\let\bss\vss\fi
  \if#1b\def\bss{\vskip\tb@bpad}\let\tss\vss\fi
  \if#1l\def\lss{\hskip\tb@lpad}\let\rss\hss\fi
  \if#1m\let\lss\hss\let\rss\hss\fi
  \if#1r\def\rss{\hskip\tb@rpad}\let\lss\hss\fi
}

\def\tb@fl{\endcell\begincell\vrule depth 0pt width \dimen0 height \dimen0 \endcell\begincell}

\@ifundefined{diagram}{}{
\def\tb@arrowpad{.5}

\newoptcommand{\tb@arrow}{\@ne}[2]{%
  \endcell
   \begingroup%
   \let\dg@getnodesize\tb@getnodesize
   \dg@USERSIZE=#1\relax%
   \ifnum\dg@USERSIZE<\@ne \dg@USERSIZE=\@ne \fi%
   \dg@parse{#2}%
   \dg@label{\tb@draw{#1}{#2}}}

\def\tb@getnodesize#1#2#3#4#5{\dimen3=\tb@arrowpad\dimen2 #4=\dimen3 #5=\dimen3\relax}
\def\tb@getnodesize#1#2#3#4#5{\ifnum#2=0\ifnum#3=0\tb@getnodesizetail{#4}{#5}\else\tb@getnodesizehead{#4}{#5}\fi\else\tb@getnodesizehead{#4}{#5}\fi}
\def\tb@getnodesizetail#1#2{\dimen3=.5\dimen2 #1=\dimen3 #2=\dimen3}
\def\tb@getnodesizehead#1#2{\dimen3=.5\dimen2 #1=\dimen3 #2=\dimen3}

\def\tb@draw#1#2#3#4{%
        \dg@X=0\dg@Y=0\dg@XGRID=1\dg@YGRID=1\unitlength=.001\dimen0%
        \dg@LBLOFF=\dgLABELOFFSET \divide\dg@LBLOFF\unitlength%
        \dg@drawcalc
        \begincell
        \let\lams@arrow\tb@lams@arrow
        \begin{picture}(0,0)\begingroup\dg@draw{#1}{#2}{#3}{#4}\end{picture}%
        \endcell
        \endgroup
        \begincell}
}

\def\tb@lams@arrow#1#2{%
 \lams@firstx\z@\lams@firsty\z@
 \lams@lastx#1\relax\lams@lasty#2\relax
 \lams@center\z@
 \N@false\E@false\H@false\V@false
 \ifdim\lams@lastx>\z@\E@true\fi
 \ifdim\lams@lastx=\z@\V@true\fi
 \ifdim\lams@lasty>\z@\N@true\fi
 \ifdim\lams@lasty=\z@\H@true\fi
 \NESW@false
 \ifN@\ifE@\NESW@true\fi\else\ifE@\else\NESW@true\fi\fi
 \ifH@\else\ifV@\else
  \lams@slope
  \ifnum\lams@tani>\lams@tanii
   \lams@ht\ten@\p@\lams@wd\ten@\p@
   \multiply\lams@wd\lams@tanii\divide\lams@wd\lams@tani
  \else
   \lams@wd\ten@\p@\lams@ht\ten@\p@
   \divide\lams@ht\lams@tanii\multiply\lams@ht\lams@tani
  \fi
 \fi\fi
 \ifH@  \lams@harrow
 \else\ifV@ \lams@varrow
 \else \lams@darrow
 \fi\fi
}

\catcode`\@=\savecatcodeat
\let\savecatcodeat\undefined

\numberwithin{equation}{section}

\makeatletter
\renewcommand{\subsubsection}{\@startsection
{subsubsection}
{3}
{0mm}
{\baselineskip}
{-0.5\baselineskip}
{\normalfont\normalsize\bfseries}}
\makeatother

\newtheorem{theorem}{Theorem}
\newtheorem{lemma}[theorem]{Lemma}
\newtheorem{proposition}[theorem]{Proposition}
\newtheorem{example}[theorem]{Example}

\newtheorem{corollary}[theorem]{Corollary}

\newtheorem{definition}[theorem]{Definition}

\newtheorem{property}[theorem]{Property}

\newtheorem{remark}[theorem]{Remark}

\def \UU {\mu/^k}
\def \cpreceq {\to_k}
\def\la{{\lambda}} 
\def\con{{ \, \subseteq \, }} 
 
\def\cal L{{\mathcal L}} 
 
\def\N{{\mathbb N}} 
\def\Z{{\mathbb Z}}

\def\ssp{\,}
\def\ses{\,=\,}
\def\ess{\:}
\def\sas{\smallskip}

\def\dg{\delta/\gamma}
\def\aa{\alpha}
\def\CY{Y}
\def \RA {\rightarrow}

\def\om {\omega}

\def\multi#1{\vbox{\baselineskip=0pt\halign{\hfil$\scriptstyle\vphantom{(_)}##$\hfil\cr#1\crcr}}}
\def\picture #1 by #2 (#3){
  \vbox to #2{
    \hrule width #1 height 0pt depth 0pt
    \vfill
    \special{picture #3} 
    }
  }

\def \core {{\mathfrak c}}
\def \kbnd { \mathfrak p}
\def \shpe { \mathfrak s}

\def\gg {\gamma}
\def\ggg {\gamma}
\def\dd {\delta}
\def\aa {\alpha}

\def \part {\vdash}
\def\La {\Lambda}
\def \scs {\ssp , \ssp}

\begin{document}

\title[Tableaux and chains in the $k$-Young lattice]
{Tableaux on $k+1$-cores, reduced words for affine permutations,
and $k$-Schur expansions}

\author{Luc Lapointe}
\thanks{Research supported in part by FONDECYT (Chile) grant \#1030114,
the Programa Formas Cuadr\'aticas of the Universidad de Talca,
and NSERC (Canada) grant \#250904}
\address{Instituto de Matem\'atica y F\'{\i}sica,
Universidad de Talca, Casilla 747, Talca, Chile}
\email{lapointe@inst-mat.utalca.cl}

\author{Jennifer Morse}
\thanks{Research supported in part by NSF grant \#DMS-0231730}
\address{Department of Mathematics, 
University of Miami, Coral Gables, Fl 33124}
\email{morsej@math.miami.edu}


\begin{abstract} 
The $k$-Young lattice $Y^k$ is a partial order on partitions with no
part larger than $k$. This weak subposet of the Young lattice originated 
\cite{[LLM]} from the study of the $k$-Schur functions $s_\lambda^{(k)}$,
symmetric functions that form a natural basis of the space spanned 
by homogeneous functions indexed by $k$-bounded partitions.  
The chains in the $k$-Young lattice 
are induced by a Pieri-type rule experimentally satisfied by the
$k$-Schur functions.  Here, using a natural bijection between $k$-bounded 
partitions and $k+1$-cores, we establish an algorithm for identifying 
chains in the $k$-Young lattice with certain tableaux on $k+1$ cores. This 
algorithm reveals that the $k$-Young lattice is isomorphic to the
weak order on the quotient of the affine symmetric group $\tilde S_{k+1}$ 
by a maximal parabolic subgroup.  From this, the conjectured $k$-Pieri rule
implies that the $k$-Kostka matrix connecting the homogeneous 
basis $\{h_\la\} _{\la\in\CY^k}$ to $\{s_\la^{(k)}\}_{\la\in\CY^k}$
may now be obtained by counting appropriate classes of tableaux on 
$k+1$-cores.  This suggests that the conjecturally positive $k$-Schur 
expansion coefficients for Macdonald polynomials 
(reducing to $q,t$-Kostka polynomials for large $k$)
could be described by a $q,t$-statistic on these tableaux, or
equivalently on reduced words for affine permutations.
\end{abstract}

\keywords{affine Weyl group, cores, $k$-Schur functions, Macdonald polynomials}

\maketitle

\section{Introduction}
\subsection{The $k$-Young lattice}
Recall that $\la$ is a successor of a partition $\mu$ in the Young lattice 
when $\la$ is obtained by adding an addable corner to $\mu$
where partitions are identified by their Ferrers diagrams,
with rows weakly decreasing from bottom to top.
This relation, which we denote  ``$\mu\RA \la$'', occurs 
naturally in the classical Pieri rule
\begin{equation}
h_1[X]\,s_\mu[X]\ses 
\sum_{\la:\, \mu\RA \la}s_\la[X]
\,,
\end{equation}
and the partial order of the Young lattice may be defined as the 
transitive closure of  $\mu\RA \la$.  It was experimentally observed that 
the $k$-Schur functions \cite{[LLM],[LM1]} satisfy the rule
\begin{equation} \label{pruleintro}
h_1[X]\,s_\mu^{(k)}[X]\ses 
\sum_{\la:\, \mu\, \RA_k\,  \la} s_\la^{(k)}[X] \, ,
\end{equation}
where ``$\mu\, \RA_k\,  \la$'' is a certain subrelation of ``$\mu\RA \la$''.
This given, the partial order of the $k$-Young lattice $\CY^k$ is defined as 
the transitive closure of $\mu\, \RA_k\,  \la$.  

The precise definition of the relation $\mu\,\RA_k\, \la$ stems from 
another ``{Schur}'' property of $k$-Schur functions.  Computational 
evidence suggests that  the usual $\omega$-involution for symmetric 
functions acts on $k$-Schur functions according to the formula
\begin{equation}
\om\,s_\mu^{(k)}[X]\ses  s_{\mu^{\om_k} }  ^{(k)}[X] \, ,
\end{equation}
where the map  $\mu\mapsto \mu^{\om_k}$ is an involution on 
$k$-bounded partitions called  ``{\it $k$-conjugation}'' that 
generalizes partition conjugation $\mu\mapsto\mu'$.  Then
viewing the covering relations on the Young lattice as 
\begin{equation}
\mu\RA\la \ess\ess\ess \Longleftrightarrow\ess\ess\ess |\la|=|\mu|+1 \quad \&
\quad \mu\, \con \, \la \quad \& \quad\mu' \, \con \, \la' \, , 
\end{equation}
we accordingly, in our previous work \cite{[LLM]},
defined $\mu\,\RA_k\,\la$ in terms of the involution 
$\mu\mapsto \mu^{\om_k}$ by 
\begin{equation} \label{eq1.3}
\mu\, \RA_k\,\la \ess\ess\ess \Longleftrightarrow\ess\ess\ess |\la|=|\mu|+1 \quad \&\quad
\mu \, \con \, \la \quad \& \quad \mu^{\om_k} \, \con \, \la^{\om_k} \, . 
\end{equation}
Thus only certain addable corners may be added to a partition $\mu$ 
to obtain its successors in the $k$-Young lattice. 
We shall call such corners the ``{\it $k$-addable corners}'' of $\mu$.

Here, we provide a direct characterization of $k$-addable corners.
This characterization is obtained by first constructing a bijection 
between $k$-bounded partitions and $k+1$-cores.  We then show that 
certain operators preserving the set of \hbox{$k+1$-cores} 
act on the $k$-Young lattice (through this bijection) by lowering or raising 
its elements according to the covering relations. Passing from an 
element to its successor by means of these operators, we thus obtain 
an algorithm for constructing any saturated chain in the $k$-Young 
lattice.  Such a construction leads us to a bijection between 
chains in the $k$-Young lattice and a new family of tableaux 
called ``$k$-tableaux" that share properties of usual semi-standard 
tableaux.

On the other hand, the weak order on the quotient of the affine symmetric 
group $\tilde S_{k+1}$ by a maximal parabolic subgroup can be characterized 
using the previously mentioned operators on cores.  This enables us to 
show that the $k$-Young lattice is isomorphic to the weak order on 
this quotient.  Consequences include a bijection between standard 
$k$-tableaux of a fixed shape and reduced words for a fixed affine 
permutation, as well as a new bijection 
between $k$-bounded partitions and affine permutations in the quotient.

\def \RRm {\Rightarrow}

To precisely summarize our results, first recall that a $k+1$-core 
is a partition with no $k+1$-hooks.  For any $k+1$-core $\ggg$,
we then define
$$
\kbnd(\gg)=(\lambda_1,\ldots,\lambda_\ell)
$$
where $\lambda_i$ is the number of cells with $k$-bounded hook length
in row $i$ of $\gamma$.  It turns out that $\kbnd(\gg)$ is a $k$-bounded
partition and that the correspondence $\gg\mapsto \kbnd(\gg)$
bijectively maps $k+1$-cores onto $k$-bounded partitions. 
With $\lambda \mapsto \core(\lambda)$ denoting the inverse of
$\kbnd$, we define the $k$-conjugation of a 
$k$-bounded partition $\la$ to be
\begin{equation} \label{newkconj}
\la^{\om_k}\ses \kbnd\big(\core(\lambda)'\big)\, .
\end{equation}
That is, if $\gamma$ is the $k+1$-core corresponding to $\lambda$, 
then $\lambda^{\omega_k}$ is the partition whose row lengths
equal the number of $k$-bounded hooks in corresponding rows of $\gamma'$.
This reveals that $k$-conjugation, which originally emerged from the action 
of the $\om$ involution on $k$-Schur functions, is none other than the 
$\kbnd$-image of ordinary conjugation of $k+1$-cores.

The $\kbnd$-bijection then leads us to a characterization for
$k$-addable corners that determine successors in the $k$-Young 
lattice.  By labeling every square $(i,j)$ in the $i^{th}$ row 
and $j^{th}$ column by its ``{\it $k+1$-residue$\, $}'', 
$j-i \mod k+1$, we find

\noindent{\bf (Theorem~\ref{khook})}
{\it Let $c$ be any addable corner of a $k$-bounded partition $\la$ 
and $c'$ (of $k+1$-residue $i$) be the addable corner of $\core(\la)$ 
in the same row as $c$.
$c$ is $k$-addable if and only if $c'$ is the highest addable corner of  
$\core(\la)$ with $k+1$-residue  $i$. 
}

\noindent
This characterization of $k$-addability leads us to a notion of 
standard $k$-tableaux which we prove are in bijection with saturated
chains in the $k$-Young lattice. 

\noindent{(\bf Definition \ref{defktab})}
{\it Let $\gg$ be a $k+1$-core and $m$ be the number of $k$-bounded 
hooks of $\gg$. A standard $k$-tableau of shape $\gg$
is a filling of the cells of $\gg$ with the letters $1,2,\ldots ,m$ 
which is strictly increasing in rows and columns and such that
the cells filled with the same letter have the same $k+1$-residue.
}

\noindent{\bf (Theorem~\ref{firstbijection})} 
{\it The  saturated chains in the $k$-Young lattice joining the 
empty partition $\emptyset$ to a given $k$-bounded partition $\la$ are 
in bijection with the standard $k$-tableaux of shape $\core(\la)$.}

We then consider the affine symmetric group $\tilde S_{k+1}$ modulo a 
maximal parabolic subgroup denoted by $S_{k+1}$.  Bruhat order on the 
minimal coset representatives of $\tilde S_{k+1}/S_{k+1}$ can be 
defined by containment of $k+1$-core diagrams (this connection is 
stated by Lascoux in \cite{[L]} and is equivalent to other 
characterizations such as in \cite{[MM],[BB]}).
From this, stronger relations among $k+1$-core diagrams can be used 
to describe the weak order on such coset representatives.  
We are thus able to prove that our new 
characterization of the $k$-Young lattice chains implies 
an isomorphism between the $k$-Young lattice and the weak order
on the minimal coset representatives.  
Consequently,  a bijection between the set of $k$-tableaux of a given shape 
$\core(\lambda)$ and the set of reduced decompositions for
a certain affine permutation $\sigma_\lambda \in \tilde S_{k+1}/S_{k+1}$ 
can be achieved by mapping:
\begin{equation}
{\mathfrak w}:
T\mapsto s_{i_\ell}\cdots s_{i_2}s_{i_1}\, ,
\end{equation}
where ${i_a}$ is the $k+1$-residue of letter $a$ in the standard 
$k$-tableau $T$.  A by-product of this result is a simple bijection 
between  $k$-bounded partitions and affine permutations in 
$\tilde S_{k+1}/S_{k+1}$:
\begin{equation} \label{mapphi}
\phi: \lambda \mapsto \sigma_\lambda \, ,
\end{equation}
where $\sigma_{\lambda}$ corresponds to the reduced decomposition
obtained by reading the $k+1$-residues of $\lambda$ from right to 
left and from top to bottom.  It is shown in \cite{[LMW]} that this 
bijection, although algorithmically distinct, is equivalent to the 
one presented by Bj\"orner and Brenti \cite{[BB]} using a notion of 
inversions on affine permutations.  It follows from our results that 
Eq.~\eqref{pruleintro} reduces simply to
\begin{equation}
h_1[X] \, s_{\phi^{-1}(\sigma)}^{(k)}[X]=
\sum_{\sigma <\!\!\! \cdot \,_w \,   \tau } s_{\phi^{-1}(\tau)}^{(k)}[X] \, ,
\end{equation}
where the sum is over all permutations that cover $\sigma$ in the
weak order on $\tilde S_{k+1}/S_{k+1}$.

As will be detailed in \S~\ref{subsym}, 
Theorem~\ref{firstbijection} also plays a role in the theory of 
Macdonald polynomials and the study of $k$-Schur functions,
thus motivating a semi-standard extension 
of Definition~\ref{defktab}:

\noindent{\bf (Definition~\ref{defktabgen})} 
{\it Let $m$ be the number of $k$-bounded hooks in a $k+1$-core $\gg$
and let $\aa=(\aa_1,\ldots,\aa_r)$ be a composition of $m$. A 
semi-standard $k$-tableau of shape $\gg$ and evaluation $\aa$ 
is a column strict filling of $\gg$ with the letters $1,2,\ldots ,r$
such that the collection of cells filled with letter $i$ is labeled with
exactly $\alpha_i$  distinct $k+1$-residues.}

As with the ordinary semi-standard tableaux, we show that there are no 
semi-standard $k$-tableaux under conditions relating to dominance order 
on the shape and evaluation.  An analogue of Theorem~\ref{firstbijection}
can then be used to show that this coincides with unitriangularity of 
coefficients in the $k$-Schur expansion of homogeneous symmetric functions 
and suggests that the $k$-tableaux should have statistics to combinatorially
describe the $k$-Schur function expansion of the Hall-Littlewood 
polynomials.  The analogue of Theorem~\ref{firstbijection} relies on the 
following 
definition: with the pair of $k$-bounded partitions $\la,\mu$
defined to be ``{\it $r$-admissible $\,$}'' if and only if 
$\la/\mu$ and  $\la^{\om_k}/\mu^{\om_k}$ are respectively horizontal
and vertical $r$-strips, we say a sequence of partitions
$$
\emptyset=\la^{(0)}\longrightarrow \la^{(1)}\longrightarrow\la^{(2)}\longrightarrow \cdots \longrightarrow \la^{(\ell)}
$$
is $\alpha$-admissible when $\la^{(i)},\la^{(i-1)}$ is a $\alpha_i$-admissible 
pair for $i=1,\ldots,\ell$.
It turns out that all $\alpha$-admissible sequences 
are in fact chains in the $k$-Young lattice and that
Theorem~\ref{firstbijection} extends to:

\noindent{\bf (Theorem~\ref{bijecgenn})}
{\it Let $m$ be the number of $k$-bounded hooks in a $k+1$-core $\gg$
and let $\aa$ be a composition of $m$. 
The collection of $\aa$-admissible chains joining $\emptyset$ to $\kbnd(\gg)$
is in bijection with the semi-standard $k$-tableaux of shape $\gg$ and 
evaluation $\aa$.
}

An affine permutation interpretation for the $\aa$-admissible 
chains that generalizes our $\mathfrak w$-bijection between standard 
$k$-tableaux and reduced words is given in \cite{[LMW]}
along with a more detailed discussion of the connection
between the type-$A$ affine Weyl group and the $k$-Schur
functions.  The reader is also referred to  \cite{[LM2]} for
a study of principal order ideals in the $k$-Young lattice 
along with further properties of the lattice such as the fact 
that the covering relation is invariant under 
translation by rectangular shapes with hook-length equal to $k$.  
This is the underlying mechanism in the proof that the $k$-Young lattice 
corresponds to a cone in a tiling of $\mathbb R^{k}$ by permutahedrons 
\cite{[Ul]}.  

As mentioned, the root of our work lies in the study of
symmetric functions.  We conclude our introduction with a 
summary of these ideas.

\subsection{Macdonald expansion coefficients} \label{subsym}

The $k$-Young lattice emerged from the experimental Pieri rule 
\eqref{pruleintro} satisfied by $k$-Schur functions.
In turn,  $k$-Schur functions have arisen from a close study of 
Macdonald polynomials. To appreciate the role of our findings in the 
theory of Macdonald polynomials we shall briefly review this connection.
To begin, we consider the Macdonald polynomial $H_\la[X;q,t]$ 
obtained from the Macdonald integral form \cite{[M]} $J_\la[X;q,t]$ 
by the plethystic substitution
\begin{equation}
H_\mu[X;q,t\, ]\ses J_\mu\big[{\textstyle{X\over 1-q}};q,t\, \big]\, .
\end{equation}
For $\mu\part n$, this  yields the Schur function expansion
\begin{equation} \label{macdoint}
H_\mu[X;q,t\, ]\ses \sum_{\la\part n} 
K_{\la \mu}(q,t) s_\la[X]\, ,
\end{equation}
where $K_{\la \mu}(q,t) \in \N[q,t] $ are known as the $q,t$-Kostka 
polynomials. Formula \eqref{macdoint}, when $q=t=1$, reduces to
\begin{equation} \label{h1nf}
h_1^n\ses \sum_{\la\part n} f_\la\, s_\la[X] \, ,
\end{equation}
where $f_\la$ is the number of standard tableaux of shape $\la$. This given,
one of the outstanding problems in algebraic combinatorics is 
to associate a pair of statistics $a_\mu(T),b_\mu(T)$ on 
standard tableaux to the partition $\mu$ so that
\begin{equation}
K_{\la \mu}(q,t)\ses \sum_{T\in ST(\la)} q^{a_\mu(T)} t^{b_\mu(T)} \, ,
\end{equation}
where ``$ST(\la)$'' denotes the collection of standard tableaux of shape $\la$.

In previous work \cite{[LLM],[LM1]}, we proposed a new approach 
to the study of the $q,t$-Kostka polynomials.  This approach is based on 
the discovery of a certain family of symmetric functions 
$\{s_\la^{(k)}[X;t]\}_{\lambda\in \CY^k}$ for each integer $k\ge 1$, 
which we have shown \cite{[LM1]} to be a basis for the space $\La^{(k)}_t$
spanned by the Macdonald polynomials $H_\mu[X;q,t\, ]$ indexed by $k$-bounded 
partitions. This revealed a mechanism underlying the structure 
of the coefficients $K_{\la\mu}(q,t)$. To be precise, 
for $\mu,\nu\in \CY^{k}$, consider 
\begin{equation} \label{mackschurintro}
H_\mu[X;q,t\, ]\ses \sum_{\nu \in \CY^{k}}  K_{\nu\mu}^{(k)}(q,t) \, 
s_\nu^{(k)}[X;t\, ]\ess\scs  \ess\text{and}
\ess\ess\ess s_\nu^{(k)}[X;t\, ]\ses \sum_\la  \pi_{\la \nu}(t)
\, s_\la[X] \,.
\end{equation}  
We then we have the factorization
\begin{equation}
K_{\la\mu}(q,t)\ses \sum_{\nu\in \CY^{k}}
\pi_{\la \nu}(t)\, K_{\nu \mu}^{(k)}(q,t)\, .
\end{equation}
It was experimentally observed (and proven for $k=2$ in \cite{[LM0],[LM1]}) 
that $K_{\nu \mu}^{(k)}(q,t) \in \N[q,t]$ and $\pi_{\la \nu}(t)\in \N[t]$.  
This suggests that  the problem of finding statistics for  $K_{\la \mu} (q,t)$ 
may be decomposed into two separate analogous problems
for $K_{\nu \mu}^{(k)}(q,t) $ and $\pi_{\la \nu}(t) $. 
We also have experimental evidence to support that 
$K_{\la \mu}(q,t)-K_{\nu \mu}^{(k)}(q,t) \in \N[q,t]$ 
which brings about the fact that $s_\la^{(k)}[X;t\,]$-expansions 
are formally simpler.

These developments prompted a close study of the polynomials  
$s_\la^{(k)}[X;1\, ] =s_\la^{(k)}[X]$. In addition to \eqref{pruleintro},
it was also conjectured that these polynomials satisfy the more general rule
\begin{equation} \label{prulesintro}
h_r[X]\,  s_\mu^{(k)}[X]\ses  
\sum_{\la/\mu=\text{horizontal $r$-strip}
\atop{ \la^{\om_k}/\mu^{\om_k}=\text{vertical $r$-strip}}}   s_\la^{(k)}[X ] 
\, .
\end{equation}
Iteration of \eqref{pruleintro} starting from  $s_{\emptyset}[X]=1$ yields 
\begin{equation} \label{h1nintro}
h_1^n[X]\ses \sum_{\la\in \CY^k} K_{\la,1^n}^{(k)} s_\la^{(k)}[X  ]\, ,
\end{equation}
while iterating \eqref{prulesintro} for suitable choices of 
$r$ gives the $k$-Schur function expansion 
of an $h$-basis element
indexed by any $k$-bounded partition $\mu$ . That is,  
\begin{equation} \label{hgenintro}
h_\mu[X]\ses \sum_{\la\in \CY^k}  K_{\la \mu}^{(k)}  \, s_\la^{(k)}[X]
\end{equation}
Since $s_\la^{(k)}[X  ]=s_\la [X  ]$ when all the hooks of $\la$ are 
$k$-bounded, we see that \eqref{h1nintro} reduces to \eqref{h1nf}
for a sufficiently large $k$. For the same reason, the coefficient 
$ K_{\la \mu}^{(k)}$ in \eqref{hgenintro} reduces to the classical 
Kostka number $ K_{\la \mu}$ when $k$ is large enough. 
Our definition of the $k$-Young lattice $\CY^k$ and admissible chains
in $\CY^k$, combined with the experimental Pieri rules \eqref{pruleintro} and 
\eqref{prulesintro}, yield the following corollary of 
Theorems~\ref{firstbijection} and \ref{bijecgenn}:

\begin{quote}
\it On the validity of \eqref{prulesintro}, the coefficient 
$ K_{\la,1^n}^{(k)}$ is equal to the number of standard $k$-tableaux 
of shape $\core(\la)$, or equivalently the number of reduced expressions 
for $\sigma_\la$, and the coefficient $ K_{\la \mu}^{(k)}$ is equal to the 
number of semi-standard $k$-tableaux of shape $\core(\la)$ and 
evaluation $\mu$.
\end{quote}
Since \eqref{mackschurintro} reduces to \eqref{h1nintro} when $q=t=1$, 
this suggests that the positivity of $K_{\la \mu}^{(k)}(q,t)$
may be accounted for by $q,t$-counting  standard 
$k$-tableaux of shape $\core(\lambda)$, or reduced words of $\sigma_{\lambda}$,
 according to a suitable statistic depending on $\mu$.  
More precisely, for $\mathcal T^k(\lambda)$ the set of $k$-tableaux
of shape $\core(\lambda)$ and $Red(\sigma)$ the reduced words for $\sigma$,
\begin{eqnarray}
H_\mu[X;q,t\, ] & \ses & \sum_{\la:\la_1\leq k} 
\left( \sum_{T\in\mathcal T^k(\lambda)} 
q^{a_\mu(T)}\,t^{b_\mu(T)}\right)
 s_\la^{(k)}[X;t\, ]\\
& \ses &
\sum_{\sigma\in\tilde S_{k+1}/S_{k+1}} 
\left( \sum_{w\in Red(\sigma)} 
q^{a_{\sigma_\mu}(w)}\,t^{b_{\sigma_\mu}(w)}\right)
 s_{\phi^{-1}(\sigma)}^{(k)}[X;t\, ] 
\, .
\end{eqnarray}

We should also mention that the relation in \eqref{hgenintro} 
was proven to be unitriangular \cite{[LM1]} with respect to 
the dominance partial order ``$\unrhd$'' as well as the
$t$-analog of this relation, given by the Hall-Littlewood
polynomials corresponding to the specialization $q=0$ of the 
Macdonald polynomials:  
\begin{equation}
H_\mu[X;0,t\, ]\ses 
\sum_{\multi{\la\in \CY^k\cr \la \unrhd \mu}} 
K_{\la \mu}^{(k)}(t) \, s_\la^{(k)}[X;t\, ] \, .
\end{equation}
Our conjecture that $K_{\la \mu}^{(k)}(q,t)\in\mathbb N[q,t]$
implies $K_{\la \mu}^{(k)}(t)$ would also have positive integer 
coefficients.  Our work here then suggests that this
positivity may be accounted for by showing that 
the $K_{\la \mu}^{(k)}(t)$ can be obtained by $t$-counting 
semi-standard $k$-tableaux according to a suitable $k$-charge 
statistic.

In conclusion, since \eqref{hgenintro} is obtained by iterating 
\eqref{prulesintro} and the resulting matrix $\|K_{\la \mu}^{(k)}\|$ 
is unitriangular, the inversion of this matrix gives a well-defined
family of functions that are conjectured to be the $k$-Schur 
functions.  This provides a relatively simple algorithm for 
computing the ``$k$-Schur functions" (at $t=1$) for anyone who 
wishes to experiment. A study of the  ``$k$-Schur functions" obtained 
in this manner is carried out in \cite{[LMW]} where 
it is shown, in particular,
that they satisfy the $k$-Pieri rule \eqref{pruleintro}.

\tableofcontents

\section{Definitions  }
\subsection{Partitions}
A partition $\lambda=(\lambda_1,\dots,\lambda_m)$ is a 
non-increasing sequence of positive integers. 
The degree of $\lambda$ is $|\lambda|=\lambda_1 +\cdots +\lambda_m$
and the length $\ell(\lambda)$ is the number of parts $m$.
Each partition $\lambda$ has an associated Ferrers diagram
with $\lambda_i$ lattice squares in the $i^{th}$ row,
from the bottom to top.  For example,
\begin{equation}
\lambda\,=\,(4,2)\,=\,
{\tiny{\tableau*[scY]{ & \cr & & & }}} \, .
\end{equation}
Given  a partition $\lambda$, its conjugate $\lambda'$ is the diagram 
obtained by reflecting  $\lambda$ about the diagonal.  A partition 
$\lambda$ is ``{\it $k$-bounded}'' if $\lambda_1 \leq k$.  Any lattice 
square in the Ferrers diagram is called a cell, where the cell $(i,j)$ 
is in the $i$th row and $j$th column of the diagram. We say 
that $\lambda \subseteq \mu$ when $\lambda_i \leq \mu_i$ for
all $i$.  The dominance order $\unrhd$ on partitions is defined by 
$\lambda\unrhd\mu$ when $\lambda_1+\cdots+\lambda_i\geq
\mu_1+\cdots+\mu_i$ for all $i$, and $|\lambda|=|\mu|$.

More generally, for $\rho \subseteq \gamma$,
the skew shape $\gg/\rho$ is identified with its
diagram $\{(i,j) : \rho_i<j\leq \gg_i\}$.
Lattice squares that do not lie in $\gg/\rho$ 
will be simply called ``{\it squares$\, $}''. We shall say that any 
$c\in \rho$ lies ``{\it below $\, $}'' $\gg/\rho$.
The ``{\it hook $\, $}'' of any lattice square $s\in \gg$ is defined as the collection of cells  of $\gg/\rho$ that lie inside the $L$ with $s$ as 
its corner. This is intended to  apply to all  $s\in \gg$ including those 
below $\gg/\rho$.  In the example below the hook of $s=(1,3)$ is 
depicted by the framed cells 
\begin{equation}
\gg/\rho\,=\,(5,5,4,1)/(4,2)\,=\,
{\tiny{\tableau*[scY]
{\cr&&\tf&\cr\bl&\bl&\tf&&\cr\bl &\bl &\bl s &\bl & \tf}}} \, .
\end{equation}
We then let $h_s(\gg/\rho)$ denote the number of cells in the hook of $s$. Thus from the example above we have
$h_{(1,3)}\big((5,5,4,1)/(4,2)\big)=3$ and 
$h_{(3,2)}\big((5,5,4,1)/(4,2)\big)=3$. 
We shall also say that the hook of a cell or a square
is $k$-bounded if its length is not larger than $k$. 

\sas

\begin{remark} \label{hooksgrowold} 
It is important to note that when row and column lengths of $\gg/\rho$ 
weakly decrease from top to bottom and left to right, then the present 
notion of hook length satisfies some of the standard inequalities of hook 
lengths.  In particular, $h_{s_1}(\gg/\rho)\geq h_{s_2}(\gg/\rho)$ 
when $s_1=(i_1,j_1)$,  
$s_2=(i_2,j_2)$ with $i_1\le i_2$
and $j_1\le j_2$ and the inequality is strict when $s_1,s_2\in\gg/\rho$ 
or $s_1\in \rho$ and $s_2\in \gg/\rho$. 
\end{remark}
Recall that a {\it ``$k+1$-core} is a partition that does not contain 
any $k+1$-hooks (see \cite{[JK]} for more on cores and residues).
The ``{\it $k+1$-residue}'' of square  $(i,j)$ is $j-i \mod k+1$.
That is, the integer in this square when squares are
periodically labeled with $0,1,\ldots,k$, where zeros lie on 
the main diagonal.  The 5-residues associated to the 5-core
$(6,4,3,1,1,1)$ are
$$
{\tiny{\tableau[scY]{\bl 4 |0|1|2,\bl 3|3,4,0,\bl 1|4,0,1,2,
\bl 3| 0,1,2,3,4,0,\bl 1}}}
$$

A cell $(i,j)$ of a partition $\gg$ with $(i+1,j+1)\not\in \gg$ is
called ``{\it extremal }''.  An extremal cell which is neither at the end 
of its row nor at the top of its column 
is called ``{\it corner extremal }''.  A ``{\it removable}'' corner of
partition $\gg$ is a cell $(i,j)\in \gg$ with $(i ,j+1),(i+1,j)\not\in \gg$
and an ``{\it  addable}'' corner of $\gg$ is a square  
$(i,j)\not\in \gg$ with $(i ,j-1),(i-1,j)\in \gg$.  All removable corners 
are extremal. We should note that $(1,\gg_1),(\ell(\gg),1)$ are removable 
corners  and  $(1,\gg_1+1),(\ell(\gg)+1,1)$ are addable. In the figure below 
we have labeled all addable corners with $a$, labeled extremal cells $e$ 
(with the corner extremals overlined), and framed the removable corners. 
\begin{equation}
{\tiny{\tableau[scY]{\bl,\bl {{a}} , \bl  ,\bl | 
\bl,e, \tf e, \bl {{a}} ,\bl |
\bl,,{\overline{e}},\tf e,\bl {{}},\bl  , \bl , \bl |
\bl ,\bl,,{{e}},\bl {{a}},\bl ,\bl |
\bl,\bl,, {\overline{e}}, e, \tf e,\bl {{a}}}}}
\end{equation}
Given any two squares, $b$ south-east of $a$, ``$a\wedge b\, $'' will 
denote the square that is simultaneously directly south of $a$ and directly 
west of $b$.

A composition $\alpha$ of an integer $m$ is a vector of positive 
integers that sum to $m$.  A ``{\it tableau}'' $T$ of shape $\lambda$ is 
a filling of $T$ with integers that is weakly increasing in rows 
and strictly increasing in columns.  The ``{\it evaluation}'' of $T$ is 
given by a composition $\alpha$ where $\alpha_i$ is the multiplicity 
of $i$ in $T$.

\def \CC {{\mathcal C}}
\def \CP {{\mathcal P}}
\def \CCk {\CC_{k+1}}

\subsection{Affine symmetric group}
The affine symmetric group $\tilde S_{k+1}$ is generated by the 
$k+1$ elements $\hat s_0,\dots, \hat s_{k}$ satisfying the affine 
Coxeter relations:
\begin{eqnarray} \label{coxeter}
\hat s_i^2  =  id, \qquad
\hat s_i \hat s_j = \hat s_j \hat s_i \quad \; (i-j\neq \pm 1\mod k+1),
\quad\text{and}\quad
\hat s_i \hat s_{i+1} \hat s_{i} =  \hat s_{i+1} \hat s_i \hat s_{i+1} \, .
\end{eqnarray}
Here, and in what follows, $\hat s_i$ is understood as
$\hat s_{i \! \! \! \mod  \! k+1}$ if $i\geq k+1$. 
Elements of $\tilde S_{k+1}$ are called affine permutations, 
or simply permutations.
A word $w=i_1 i_2\cdots i_m$ in the alphabet $\{0,1,\ldots,k\}$
corresponds to the permutation $\sigma \in \tilde S_{k+1}$ if
$\sigma=\hat s_{i_1} \dots \hat s_{i_m}$.  The {\it ``length"} of 
$\sigma$, denoted $\ell(\sigma)$, is the length of the shortest 
word corresponding to $\sigma$.  Any word for $\sigma$ with
$\ell(\sigma)$ letters is said to be {\it ``reduced"}. 
We denote by $Red(\sigma)$ the set of all reduced words of $\sigma$.

The weak order on $\tilde S_{k+1}$ is defined through the following 
covering relations:
\begin{equation}
\sigma <\!\!\!\!\cdot_{w} \tau  \quad \Longleftrightarrow \quad  \tau = \hat s_i \, \sigma 
\text{ for some } i\in \{0,\dots,k\}\, , \text{ and } \ell(\tau)>\ell(\sigma) \, ,
\end{equation}
while the Bruhat order is such that $\sigma<_b\tau$ if there exist reduced words
$w$ and $u$, corresponding to $\sigma$ and $\tau$ respectively, such that $w$ is
a subword of $u$.

The subgroup of $\tilde S_{k+1}$ generated by the subset 
$\{\hat s_1,\ldots,\hat s_{k}\}$ is a maximal parabolic subgroup
that is isomorphic to the symmetric group.   We thus denote this
subgroup by $S_{k+1}$ and shall consider the set of minimal coset 
representatives of $\tilde S_{k+1}/S_{k+1}$.  It is important to
note that if $\sigma$ is not the identity, then
$\sigma\in \tilde S_{k+1}/S_{k+1}$ if and only if every $w \in Red(\sigma)$ 
ends in a zero.  
That is, the reduced expressions are all of the form
$w={i_1} \dots {i_{m-1}}\, {0}$.

\section{Bijection: $k+1$-cores and $k$-bounded partitions}

Let $\CC_{k+1}$ and $\CP_k$ respectively denote the collections of $k+1$ 
cores and $k$-bounded partitions. We start by showing that a bijection
between these sets can be defined by the map
$$
\kbnd : \gg \to (\lambda_1,\ldots,\lambda_\ell)
$$
where $\lambda_i$ is the number of cells with a $k$-bounded hook
in row $i$ of $\gamma$. If $\rho(\gg)$ is the partition consisting only of the 
cells in $\gg$ whose hook lengths exceed $k$, then $\kbnd(\gg)=\lambda$ is
equivalently defined by letting $\lambda_i$ denote the length of 
row $i$ in the skew diagram $\gg/\rho(\gg)$. For example, with $k=4$:
\begin{equation*}
\gamma = 
{\tiny{\tableau*[scY]{  
\cr  
\cr  & 
\cr  &  & 
\cr  & & & & 
\cr &&  & & & & & &
}}}
\qquad
\gamma/\rho(\gg)= {\tiny{\tableau*[scY]{  
\cr  
\cr  & 
\cr \bl &  & 
\cr  \bl &\bl & & & 
\cr \bl &\bl& \bl & \bl &\bl & & & &
}}}
\qquad
\qquad
{\kbnd(\gg)}
= {\tiny{\tableau*[scY]{  
\cr  
\cr  & 
\cr    & 
\cr   & & 
\cr   & & &
}}}
\end{equation*}
Although it is not immediate that the codomain of $\kbnd$ is 
$\mathcal P_k$, we shall find that each diagram $\gamma/\rho(\gamma)$
can be uniquely associated to a skew diagram constructed from some
$k$-bounded partition $\lambda$.

\begin{definition} \label{def5}
For any $\lambda\in\mathcal P_k$, the
``$k$-skew diagram of $\lambda$" is the diagram $\lambda/^k$
where

(i) row $i$ has length $\lambda_i$ for $i=1,\ldots,\ell(\lambda)$

(ii) no cell of $\lambda/^k$ has hook-length exceeding $k$

(iii) all squares below $\lambda/^k$ have hook-length exceeding $k$.

\noindent
\end{definition}
We shall thus find that $\kbnd$ is a bijection from $\mathcal C_{k+1}$ to
$\mathcal P_k$ with inverse $\core$ defined:
\begin{definition}
For $\lambda$ a $k$-bounded partition and
$\lambda/^k=\gamma/\rho$, define $\core(\lambda)=\gamma$.
\end{definition}

To this end, we start by characterizing the skew diagrams 
$\gg/\rho(\gg)$ by the following lemma, and consequently 
find that $\kbnd(\gamma)\in\mathcal P_k$.

\begin{lemma} \label{propo3.1}
Let $\rho \, \con \,  \gg$ be partitions.  Then 
$\gg$ is a $k+1$-core and $\rho=\rho(\gg)$ if and only if the
skew partition $\gg/\rho$ has the following four properties:
\begin{itemize}
\item[{\rm (i)}] the row lengths   of $\gg/\rho$ weakly decrease from bottom to top,

\item[{\rm (ii)}] the  column lengths of $\gg/\rho$ weakly decrease  
from  left to right,

\item[{\rm(iii)}] the hooks of the cells of $\gg/\rho$ have at most $k$ cells,

\item[{\rm(iv)}] the  squares below $\gg/\rho$ have hook-lengths exceeding $k$. 
\end{itemize}
\end{lemma}

\begin{proof}
We shall prove that conditions (iii) and (iv) are sufficient.
Let $c\in\rho$.  Condition (iv) asserts that the hook of $c$ 
in $\gg/\rho$ contains at least $k+1$ cells.  
Since the hook of $c$ in $\gg$ contains at least these cells
and $c$ itself, $h_c(\gg)>k+1$.  Moreover (iii) assures that 
all cells of $\gg/\rho$ have hook length $\le k$. Thus $\gg$ has no $k+1$ 
hooks and therefore is a $k+1$-core. Since all the cells 
of $\gg$ whose hook length exceeds $k$ lie in $\rho$, it follows 
that $\rho=\rho(\gg)$ as desired.

To show that conditions (i)-(iv) are necessary, let $\gg$ be a $k+1$-core 
and $\rho=\rho(\gg)$.  Refer below to a typical case of two successive 
rows in $\gg$, where $b$ is the cell at the end of row $i+1$,
$c$ is the cell at the end of row $i$, and $a$ labels the cell at the top of 
column $\rho_i$ in $\gg$. Thus $a\wedge c$ is the cell $(i,\rho_i)$ 
in $\rho$ (labeled with a ``$\bullet$''). 
Let $d$ be any cell of row $i+1$ that has at least
$\gg_i-\rho_i$ cells of $\gg$ to its right.
\begin{equation}
{\tiny{\tableau*[scY]
{ \bl & \bl & \bl & a \cr\bl  & \bl  & \bl  & \cr \bl & \bl & \bl& \cr
\bl &\bl &\bl &\cr d & & & & &  &  &  &  b\cr
\bl & \bl  & \bl  & \bl \bullet   & & & & & & & & c\cr
}}} 
\end{equation}
The definition of $\rho(\gg)$ implies that
$h_{a\wedge c}(\gg)\ge k+2$ since $\gg$ is a $k+1$-core. 

To show (i) we must prove that $\gg_i-\rho_i\ge \gg_{i+1}-\rho_{i+1}$,
or equivalently that all such $d's$ are in $\rho$. 
For this, we observe that the hook length of $d$ in $\gg$
is at least equal to the hook length of $a\wedge c$ minus one.
Thus, $h_{a\wedge c}(\gg)\geq k+2$ implies that the hook length 
of $d$ is at least $k+1$ and therefore $d$ belongs to $\rho$.
Next note that since the conjugate of a $k+1$ core is also a $k+1$ core 
and $\rho(\gg')=\rho(\gg)'$, condition (ii) for $\gg/\rho(\gg)$ follows 
from (i) for $\gg'/\rho(\gg')$.

Condition (iii) is an immediate consequence of the definition of $\rho(\gg)$. 
To prove (iv), consider $a\wedge c$ (the ``$\bullet$'' in our figure) 
as a typical removable corner of $\rho(\gg)$. Since 
$h_{a\wedge c}(\gg)\geq k+2$, at least $k+1$ of these cells
lie in $\gg/\rho(\gg)$ implying $h_{a\wedge c}(\gg/\rho)\geq k+1$.
Since every square below  $\gg/\rho(\gg)$ is weakly {\it south-west} 
of such a removable corner in $\rho(\gg)$, condition (iv) follows 
Remark~\ref{hooksgrowold} given (i) and (ii).
\end{proof}

We thus have that $\kbnd$ maps   $\CCk$ into  $\CP_k$ since the parts 
of $\kbnd(\gg)$ are weakly increasing by condition (i) and do not 
exceed $k$ by condition (iii).  To show that this map is a bijection, 
we will identify its inverse by considering the following auxiliary 
result:
\begin{lemma} \label{propo3.2}
For any $k$-bounded partition $\la=(\la_1,\la_2,\ldots,\la_r)$,
there is a unique sequence of skew diagrams
$\lambda^{(r)}\!\!/^k,\lambda^{(r-1)}\!\!/^k,\ldots,\lambda^{(1)}\!\!/^k$ 
where $\lambda^{(r)}\!\!/^k=(\la_r)$  and $\lambda^{(i)}\!\!/^k$ 
is obtained by attaching a row of length $\la_i$ to
the bottom of $\lambda^{(i+1)}\!\!/^k$ such that:
\begin{itemize}
\item[{\rm (1)}]
the hook lengths of $\lambda^{(i)}\!\!/^k$ do not exceed $k$ 
\item[{\rm (2)}] all the lattice squares below $\lambda^{(i)}\!\!/^k$ 
have hook lengths exceeding $k$. 
\end{itemize}
In particular, $\la/^k=\lambda^{(1)}\!\!/^k$ 
is the unique skew partition $\gg/\rho$ 
such that 
\begin{itemize}
\item[{\rm (a)}] the row lengths of $\gg/\rho$ are the parts of $\la$

\item[{\rm (b)}] $\gg$ is a $k+1$-core and $\rho=\rho(\gg)$
\end{itemize}
\end{lemma}
\def \ogg {\overline \gg}
\def \odd {\overline \rho}

\begin{proof}
To prove that conditions (1) and (2) uniquely determine $\lambda^{(i)}\!\!/^k$ 
from $\la^{(i+1)}\!\!/^k$, let \hbox{$\la^{(i+1)}\!\!/^k=\gg/\rho$} 
with $\gg=(\gg_{i+1},\gg_{i+2},\ldots ,\ggg_r)$ and 
$\rho=(\rho_{i+1},\rho_{i+2},\ldots ,\rho_r)$. Inductively assume 
that all conditions have been met up to this point. 
By construction, we have $\la^{(i)}\!\!/^k=\ogg/\odd$ with 
\begin{equation} \label{eq3.1}
\ogg=(a +\la_i,\gg_{i+1},\gg_{i+2},\ldots ,\ggg_r)\ess\ess\hbox{and}\ess\ess 
\odd=(a,\rho_{i+1},\rho_{i+2},\ldots ,\rho_r) 
\end{equation}
for some $a\geq \rho_{i+1}$.  We claim that conditions (1) and (2) 
uniquely determine $a$.  From \eqref{eq3.1} we derive that the 
hook length of the first cell in the bottom row of $\la^{(i)}\!\!/^k$ 
(the leftmost framed cell in the figure) is $b_{a+1}+\la_i$
where $b_j$ is the length of the $j^{th}$ column of $\la^{(i+1)}\!\!/^k$. 
\begin{equation}
{\tiny{\tableau*[scY]
{  \cr    \cr   &   &   \cr  \fl & 
& & & \bl & \bl &\bl {a+1}  \cr
\fl & & & & \bl & \bl & \bl  \downarrow \cr \fl & \fl& \fl &  &  &    \cr
\fl & \fl& \fl & \fl &  &  &  &  \cr
 \bullet  & \bullet  & \bullet   & 
\bullet & X & X &\tf &\tf &\tf  &\tf  
& \bl & \bl \leftarrow & \bl i 
}}} 
\end{equation}
To satisfy (1) we must have 
\begin{equation} \label{eq3.2}
b_{s}+\la_i\le k\, ,\quad \quad {\rm for~all~} s \geq a+1 \, .
\end{equation}
To satisfy (2), the squares {\it west} of the added row 
must have hook lengths $\ge k+1$.  That is,
\begin{equation} \label{eq3.3}
b_{s}+\la_i\ge k+1\, ,\quad \quad {\rm for~all~} s \leq a \,  .
\end{equation}
Since $b_{a}+\la_i\ge b_{a}+\la_{i+1}$, the inductive hypothesis 
guarantees that \eqref{eq3.3} will be true for all $a\leq \rho_{i+1}$. 
That is, the squares marked with a ``$\bullet$'' in the figure 
will necessarily have hook lengths exceeding $k$.
It follows from these observations that to obtain 
a skew shape that satisfies both (1) and (2), we are forced to 
take $a$ as the smallest integer such that $b_{a}+\la_i\le k$.
In this case, \eqref{eq3.3} is automatically satisfied.  
And \eqref{eq3.2} follows 
because $b_s \geq b_{a+1}$ for all $s \geq a+1$,
given that when considering only the columns of $\la^{(i+1)}\!\!/^k$ 
starting from column $a+1$, the diagram $\la^{(i+1)}\!\!/^k$
is that of a partition.  This completes the induction.

Now let $\la^{(1)}\!\!/^k=\gg/\rho$. Our construction assures property (a). 
Also by construction, $\gg/\rho$ satisfies conditions (iii) and (iv) 
of Lemma~\ref{propo3.1}, which were shown in the proof
of that lemma to be sufficient to guarantee that 
$\gg$ is a $k+1$-core with $\rho=\rho(\gg)$.  Conversely,
if a skew diagram $\ogg/\odd$ satisfies (a) and (b)
then (a) implies $\ogg=(\ogg_1,\ogg_2,\ldots \ogg_r)$ and 
$\odd=(\odd_1,\odd_2,\ldots \odd_r)$ with $\ogg_i-\odd_i=\la_i$.
Moreover property (b) assures that all the hook lengths of 
$\ogg/\odd$ do not exceed $k$ and all the squares below $\rho$
have lengths exceeding $k$. These two properties thus necessarily 
hold for all the skew diagrams 
$$
{\overline \la}^{(i)}\!\!/^k\ses 
(\ogg_i,\ogg_{i+1},\ldots \ogg_r)/(\odd_i,\odd_{i+1},\ldots \odd_r)
\, .
$$
Therefore, $\ogg/\odd={\la}^{(1)}\!\!/^k$
by the uniqueness of such diagrams satisfying (1)-(2).
\end{proof}

Note that the proof of Lemma~\ref{propo3.2} reveals that 
$\lambda^{(i)}\!\!/^k$ is the diagram obtained by attaching the 
row of length $\lambda_i$ to the bottom of $\lambda^{(i+1)}\!\!/^k$ 
in the leftmost position so that no hook-lengths exceeding $k$ are 
created. This algorithm gives a convenient method for constructing
$\lambda/^k$.

\begin{example} \label{exskew}
Given $\lambda =(4,3,2,2,1,1)$ and $k=4$,
\begin{equation*}
\lambda = {\tiny{\tableau*[scY]{  \cr  \cr  & \cr &
\cr & & \cr & & & }}}
\quad
\implies
\quad
\lambda/^4 = {\tiny{\tableau*[scY]{  
\cr  
\cr  & 
\cr \bl &  & 
\cr  \bl &\bl & & & 
\cr \bl &\bl& \bl & \bl &\bl & & & &
}}}
\qquad\text{revealing that}
\;
\core(\lambda)= {\tiny{\tableau*[scY]{  
\cr  
\cr  & 
\cr  &  & 
\cr  & & & & 
\cr  && & & & & & &
}}}
\end{equation*}
\end{example}

\medskip

We are now in the position to prove our bijection:

\begin{theorem} \label{theo6moi}
$\kbnd$ is a bijection from $\CCk$ onto $\CP_k$ with 
inverse $\core$.
\end{theorem}
\begin{proof} 
For $\lambda\in\mathcal P_k$, consider $\gg=\core(\lambda) \in {\mathcal C}_{k+1}$.
Since $\lambda/^k=\gg/\rho$ by definition of $\core$,
Lemma~\ref{propo3.2}(a) and (b) implies that
$\kbnd(\gg)=\lambda$ and thus $\kbnd\left(\core(\lambda)\right)=\lambda$.  
Now consider $\gg\in\mathcal C_{k+1}$.  Lemma~\ref{propo3.1} 
implies that $\kbnd(\gg)=(\gg_1-\rho(\gg)_1,\ldots,\gg_n-\rho(\gg)_n)
\in\mathcal P_k$, and thus by Lemma~\ref{propo3.2},
\begin{equation}
\label{b2b3}
\gg/\rho(\gg)=\kbnd(\gg)/^k
\,.
\end{equation}
Therefore, $\core\left(\kbnd(\gg)\right)=\gg$
by definition of $\core$.
\end{proof}

\section{The $k$-lattice}

The notion of a $k$-skew diagram gives rise to an involution on 
$\mathcal P_k$, where $\lambda$ is sent to the partition whose rows 
are obtained from the columns of $\lambda/^k$:
\begin{definition}
For any $\lambda\in\mathcal P_k$, the ``{\it $k$-conjugate}'' 
of $\lambda$ denoted $\lambda^{\omega_k}$ is the $k$-bounded 
partition given by the columns of $\lambda/^k$.
\label{defconju}
\end{definition}
Equivalently, we may define the $k$-conjugation as the partition given 
by the number of $k$-bounded cells in the columns of $\core(\lambda)$,
or simply $\lambda^{\omega_k}=\kbnd(\core(\lambda)')$.  This given, 
since $\kbnd=\core^{-1}$, we see that $k$-conjugation is 
an involution by:
\begin{equation}
(\lambda^{\omega_k})^{\omega_k} = \kbnd\Bigr(\left[ 
\core\bigl(\kbnd\bigr(\core(\lambda)' \bigl)\bigr)\right]'\Bigr)
= \kbnd\bigl(\left[\core(\lambda)' \right]' \bigl)= 
\kbnd\bigl(\core(\lambda) \bigr)=\lambda \, .
\end{equation}

\begin{example}
With $\lambda$ as in Example~\ref{exskew},
the columns of $\lambda/^4$ give
$\lambda^{\omega_4} = (3,2,2,1,1,1,1,1,1)$.
\end{example}

\begin{remark}
If $h_{(1,1)}(\lambda)\leq k$, all hooks of $\lambda$
are $k$-bounded and thus $\lambda/^k=\lambda$.
In this case, $\lambda^{\omega_k}=\lambda'$.
\label{re2}
\end{remark}

Now, we can consider a partial order 
``$\, \preceq\, $'' on the collection of $k$-bounded 
partitions stemming from $k$-conjugation.

\begin{figure}[htb]
\epsfig{file=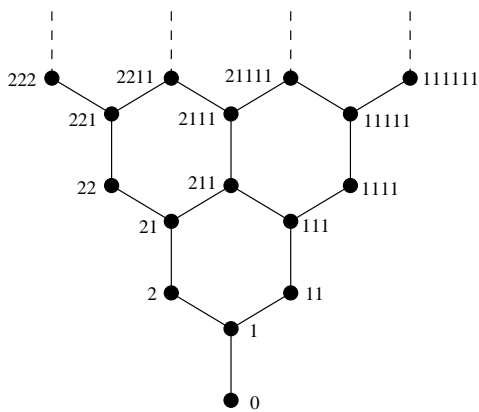}
\caption{Hasse diagram of the $k$-Young lattice in the case $k=2$.}
\label{posetexample}
\end{figure}

\begin{definition}
\label{korder}
The ``{\it $k$-Young lattice}'' $\preceq$  on partitions in 
$\mathcal P_k$ is defined by the covering relation
\begin{equation}
\lambda\cpreceq\mu\quad\text{when}
\quad
\lambda \subseteq \mu \quad\text{ and}\quad
\lambda^{\omega_k}\subseteq \mu^{\omega_k}
\end{equation}
for $\mu,\lambda\in\mathcal P_k$
where $|\mu|-|\lambda|=1$. 
Figure~\ref{posetexample} gives the case $k=2$.
\end{definition}

While this poset on $k$-bounded partitions originally arose in 
connection to a rule for multiplying generalized Schur functions 
\cite{[LLM]}, we shall see in \S~\ref{weak} that this poset 
turns out to be isomorphic to the weak order on the quotient 
of the affine symmetric group by a maximal parabolic subgroup.  
Consequently, the $k$-Young lattice is in fact a lattice 
\cite{[W]} (see \cite{[Ul]} for a proof that follows from the identification
of the $k$-Young lattice as a cone in a tiling of $\mathbb R^k$ by
permutahedrons).  

The $k$-Young lattice generalizes the Young lattice since:

\begin{property}
\label{kord}
$\lambda\preceq\mu$ reduces to $\lambda\leq \mu$
when $\lambda$ and $\mu$ are partitions such that 
$h_{(1,1)}(\mu)\leq k$.
\end{property}

\begin{proof}
Since $\lambda\subseteq\mu$ when $\lambda\preceq\mu$,
$h_{(1,1)}(\mu)\leq k$ implies that
$h_{(1,1)}(\lambda)\leq k$.  Remark~\ref{re2} then implies that
$\lambda^{\omega_k}=\lambda'$ and $\mu^{\omega_k}=\mu'$.
Thus, the conditions that $\lambda \subseteq \mu$ and
$\lambda^{\omega_k}\subseteq\mu^{\omega_k}$
reduce to $\lambda \subseteq \mu$ and $\lambda' \subseteq \mu'$,
or simply $\lambda\subseteq \mu$.
\end{proof}

Although the ordering $\preceq$ is defined by the covering relation
$\cpreceq$, the definition implies that

\begin{property}
\label{contains}
If $\lambda \preceq \mu$, then $\lambda \subseteq \mu$ and
$\lambda^{\omega_k}\subseteq \mu^{\omega_k}$.
\end{property}

It is important to notice that the converse of this statement 
does not hold.  For example, with $k=3$, $\lambda=(2,2)$ and
$\mu=(3,2,1,1,1,1)$, we have $\lambda^{\omega_k}=\lambda$ and 
$\mu^{\omega_k}=\mu$ satisfying $\lambda \subseteq \mu$ and
$\lambda^{\omega_k}\subseteq \mu^{\omega_k}$, but 
$\lambda \not \preceq \mu$.  This is easily seen from 
\cite[Th. 19]{[LM2]} (since $\lambda$ contains the 3-rectangle
$(2,2)$ while $\mu$ does not), or can be more tediously verified 
by constructing all chains using Theorem~\ref{khook}.

Similarly, although ${\core}(\la)\con {\core}(\mu)$ does not 
necessarily imply $\la\con \mu$, the converse is nevertheless true. 

\begin{property} \label{onedirect}
$\la\con\mu$ implies ${\core}(\la)\con {\core}(\mu)$. 
\end{property}
\begin{proof}
Let $\gg_\la={\core}(\la)$ and $\gg_\mu= {\core}(\mu)$.  Assume by 
contradiction that there is some row of $\gg_\mu$ that is strictly 
shorter than a row of $\gg_\la$ and let row $r$ be the highest such 
row. Let $s$ be the last cell of row $r$ of $\rho(\gg_\la)$ and $s'$ be 
the first cell of row $r$ of $\gg_\mu/\rho(\gg_\mu)$. Since 
$\gg_\la$ has $\la_r$  $k$-bounded hooks in row $r$ and $\gg_\mu$ has 
$\mu_r\ge\la_r$, the choice of $r$ forces $s'$ to be weakly to the 
left of $s$.  Now let $\ell$ be the number of cells  above $s$ in $\gg_\la$
and $\ell'$ be the number of cells above $s'$ in $\gg_\mu$. 
Since all rows of $\gg_\la$ above row $r$ are weakly smaller than the 
corresponding rows of $\gg_\mu$ and  $s'$ is weakly to the left of $s$,
we must have $\ell'\ge \ell$. Thus the choice of $s$ and $s'$ gives 
$$
k\ge h_{s'}(\gg_\mu)=\ell'+\mu_r\ge \ell+\la_r=h_{s}(\gg_\la)-1>k\,,
$$
where the last inequality holds since $h_{s}(\gg_\la)>k$ and
$\gg_\la$ has no $k+1$-hooks.  The result follows by
contradiction.
\end{proof}

In what follows, we shall develop an explicit description 
of the chains in this poset and provide a bijection with certain
tableaux.  These tableaux will then play a central role in the 
connection between the $k$-Young lattice and the weak order,
and will also be discussed in our study of Macdonald polynomials 
(see \S~\ref{symfu}).

\section{$k+1$-cores}

Since the set of $\mu$ such that $\mu\supset\lambda$ and 
$|\mu|=|\lambda|+1$ consists of all partitions obtained by adding 
a corner to $\lambda$, a subset of these partitions will be
the elements that cover $\lambda$ with respect to $\preceq$.
The definition of $\preceq$ implies that to determine which corners 
can be added to give partitions that cover $\lambda$, we must find which 
corners can be added to $\lambda$ so that the resulting diagram has a 
$k$-conjugate diagram that differs from $\lambda^{\omega_k}$ by only one box.  
Since a $k$-conjugate diagram is given by the number of $k$-bounded 
cells in the columns of a $k+1$-core, a close study of $k+1$-cores will 
enable us to characterize the allowable, or ``{\it $k$-addable}", corners.

We begin with a number of basic properties of cores that rely
on their associated residues.  For the sake of completeness,
we include all proofs although some may be known.
For any integer $d$, we shall consider the diagonals 
of a partition, $D_d\ses \{(i,j)\, :\, j-i=d \, \}$.
Note that a fixed $k+1$-residue $r=0,1,\ldots k$ occurs in 
successive diagonals  $D_{r+\ell(k+1)}$ for any integer $\ell$.
A sequence of lattice cells $c_0,c_1, \ldots, c_n$
forms a ``{\it $k+1$-string }''  if the cells  respectively lie 
in the successive diagonals:
$$
D_{r+i(k+1)}, D_{r+(i+1)(k+1)}, D_{r+(i+2)(k+1)},\ldots,D_{r+(i+n)(k+1)}\, .
$$
As such, all cells in a $k+1$-string have the same $k+1$-residue.
For any $0\leq r\leq k$, we shall also say that a square $c'$ on the 
diagonal $D_{r+(\ell)(k+1)}$ is a ``$k+1${\it -predecessor}'' of any 
square on $D_{r+(\ell+1)(k+1)}$.  It is important to notice that
if cell $c'$ is a $k+1$-predecessor of a cell $c$ in a 
partition $\lambda$, then $h_{c'\wedge c}(\lambda)=k+2$.   
This given, a $k+1$-string of cells $c_0,c_1,\ldots,c_n$ is 
simply a succession of cells 
where $c_i$ is a $k+1$-predecessor of $c_{i+1}$.  

Our point of departure here is the 
following known \cite{[GKS]} basic result.

\begin{property} \label{propo4.1}
Let $\gg$ be a $k+1$-core.
\begin{itemize}
\item[{\rm (1)}] 
Let $c$ and $c'$ be extremal cells of $\gg$ with the same $k+1$-residue
($c'$ weakly north-west of $c$).  
\begin{itemize}
\item[{(a)}] If $c$ is at the end of its row, then so is $c'$.
\item[{(b)}] If $c$ has a cell above it, then so does $c'$.
\end{itemize}
\item[{\rm (2)}] Let $c$ and $c'$ be extremal cells of 
$\gg$ with the same $k+1$-residue ($c'$ weakly south-east of $c$).  
\begin{itemize}
\item[{(a)}] If  $c$ is at the top of its column, then so is $c'$.
\item[{(b)}] If $c$ has a cell to its right, then so does $c'$.  
\end{itemize}
\item[{\rm (3)}] Let $c$ be a corner extremal cell and
$c'$ be an extremal cell of the same $k+1$-residue as $c$.
\begin{itemize}
\item[{(a)}] 
If $c'$ is weakly south-east of $c$, then $c'$ has a cell to its right.
\item[{(b)}] If $c'$ is weakly north-west of $c$,
then $c'$ has a cell above it.
\end{itemize}
\end{itemize}
\end{property}
\begin{proof} 
1(a): Given that there is no cell to the right of $c$, it 
suffices to prove that there is no cell $x$ to the right 
of the extremal cell $c'$ that is a $k+1$-predecessor of $c$ -- 
by iteration the property will follow for non-predecessors $c'$.  
If $x\in\gamma$ then the hook-length 
of the cell determined by $x$ and $c$ is $k+1$ since no cell
lies above $x$ (it is to the right of an extremal cell).
However, this contradicts that $\gamma$ is a $k+1$-core
implying that there is no cell to the right of $c'$. 
$$
{\tiny{\tableau[scY]{,,c',\bl x|,\bullet,,,,|\bullet,,,,,c|
,,,,\bullet,,, |
}}} 
$$
2(a) follows from 1(a) since the transpose of a $k+1$-core is a $k+1$-core.
Further, 1(b) and 2(b) are simply the contrapositive of 2(a) and 1(a) 
respectively with 
$c \leftrightarrow c'$.  Finally, since a corner extremal cell 
has a cell to its right and above it, 3(a)
and 3(b) follow respectively from 2(b) and 1(b).
\end{proof}

\begin{remark}
\label{both}
A $k+1$-core $\gg$ never has both a removable corner and an addable corner
of the same $k+1$-residue.  This follows by assuming there is an 
addable corner $c$ of some $k+1$-residue $i$ and using 
Property~\ref{propo4.1}(3)
with the corner extremal cell $e$ immediately south-west of $c$. 
The proposition gives
that all extremals of 
$k+1$-residue $i$ either have a cell to their right or above.  
Therefore they are not removable corners.
\end{remark}

\def \oc {\overline o}

\begin{property} \label{propodisjoint}
Let $\gg$ be a $k+1$-core
\begin{itemize}
\item[{(i)}]
If $\gg$ has a removable corner of $k+1$-residue $i$,
then the collection of all removable corners of $\gg$ with 
$k+1$-residue $i$ forms a $k+1$-string.
\item[{(ii)}] If $\gg$ has an addable  corner with $k+1$-residue $i$,
then the collection of all addable corners of $\gg$ with $k+1$-residue $i$ 
forms a $k+1$-string.
\end{itemize}
\end{property}
\begin{proof}
Let $c$ and $c'$ be the leftmost and rightmost removable corners 
of $\gg$ with $k+1$-residue $i$ and let 
$c=c_0\scs c_1 \scs c_2 \scs \ldots \scs c_m=c'$
be the extremal cells (with $c_{j}$ a \hbox{$k+1$-predecessor} 
of $c_{j+1}$) lying between $c$ and $c'$.  
By Property~\ref{propo4.1}(2a), each $c_j$ lies at the top its column 
since it is south-east of $c_0$ and by (1a) lies at the end of its row 
since it is north-west of $c_m$.  Therefore, each $c_j$ 
is a removable corner.

To prove (ii), now let $c$ and $c'$ be leftmost and rightmost addable corners 
with $k+1$-residue $i$, and let $e$ and $e'$ be the extremal cells
immediately south-east of $c$ and $c'$ respectively. With
$e=e_0\scs e_1 \scs e_2 \scs \ldots \scs e_m=e'$ the extremal cells 
(with $e_{j}$ a \hbox{$k+1$-predecessor} of $e_{j+1}$) lying between 
$e$ and $e'$, we claim that each $e_j$ is corner extremal.  
This follows from Property~\ref{propo4.1}(3). Indeed, 
each $e_j$ must have a cell to its right because it is south-east of 
$e$ and must have a cell above because it is north-west of $e'$.  
This forces the square $c_j$ that is immediately north-east of $e_j$ 
to be an addable corner of $\gg$. Since each $c_j$ has $k+1$-residue $i$ and 
is a $k+1$-predecessor of $c_{j+1}$ it follows that
$c=c_0\scs\ c_1 \scs  c_2 \scs \ldots \scs   c_{m}= c'$
forms a $k+1$-string with head $c$ and tail $c'$. 
\end{proof}

Armed with these special properties of $k+1$-cores, we turn 
to the study of certain operators that help us characterize 
the $k$-addable corners in the $k$-order, and that
enable us to identify the $k$-Young lattice with the 
weak order on $\tilde S_{k+1}/S_{k+1}$.  Operators that 
add a corner of given residue to partitions arose in 
\cite{[DJKMO]} and \cite{[MM]},
and coincide with those introduced in \cite{[Ro]}.  
In the case of a $k+1$-core, since
there is never both a removable and 
addable corner with the same $k+1$-residue by 
Remark~\ref{both}, we consider the operator \cite{[L]}
that deletes or adds all such corners from 
elements in $\mathcal C_{k+1}$.  That is,

\begin{definition}
\label{kcore}
The {\it ``operator $s_i$"} acts on a $k+1$-core by:
\begin{itemize}
\item[{(a)}] removing all removable corners with $k+1$-residue $i$ if
there is at least one removable corner of $k+1$-residue $i$

\item[{(b)}] adding all addable corners with $k+1$-residue $i$ if
there is at least one addable corner with $k+1$-residue $i$ 

\item[{(c)}] leaving it invariant when there are no addable or removable corners
of $k+1$-residue $i$.
\end{itemize}
\end{definition}

We now give a number of properties that concern the $s_i$ operators,
beginning with the observation that they preserve the
set $\mathcal C_{k+1}$.  Note, some properties given here
are implied in \cite{[L]}, but we shall
include all proofs for the sake of completeness.

\begin{property} 
\label{propouter}
Let $\gamma$ be a $k+1$-core.  
\begin{itemize}
\item[{(i)}]
If $\gamma$ has an addable corner of $k+1$-residue $i$, 
then $s_i(\gamma)$ is a $k+1$-core whose
shape is obtained by adding all addable corners of $k+1$-residue
$i$ to $\gamma$.  
\item[{(ii)}]
If $\gamma$ has a removable corner of $k+1$-residue $i$, then 
$s_i(\gamma)$ is a $k+1$-core
obtained by deleting all removable corners of $k+1$-residue $i$
from $\gamma$.
\end{itemize}
\end{property}
\begin{proof}
Let $c_1,\dots,c_n$ denote the collection of addable corners 
of $\gamma$ with $k+1$-residue $i$ where $c_j$ is 
a $k+1$-predecessor of $c_{j+1}$ for $j=1,\dots,n-1$.  
By Definition~\ref{kcore}, the diagram of $s_i(\gamma)$ 
is obtained by adding $c_1,\ldots,c_n$ to $\gamma$.  
Since $\gamma$ is a $k+1$ core and no hook of $\gamma$ is affected 
by the action of $s_i$ unless it corresponds to a cell that lies in
a column or row containing some $c_j$,  it suffices 
to check that there are no $k+1$-hooks in the rows 
of $s_i(\gamma)$ containing $c_j$ (the columns will
have no $k+1$-hooks by the transpose of our argument).
$$
{\tiny{\tableau[scY]{ , \bl b ,\bl a| ,\bullet, |\bullet,,,|,,x,,,,,\bl c_1|
,,,,,,\bullet,|,,,,,\bullet,,,,}}}
$$
First, there can only be a $k+1$-hook in the row of $s_i(\gamma)$ 
containing $c_1$ if $\gamma$ has a $k$-hook in this row. 
Assume by contradiction that a cell $x$ in the row 
with $c_1$ has a $k$-hook in $\gamma$ (see the figure). 
Let $a\not\in\gamma$ denote the lowest square 
at the top of the column containing $x$.
Since $a$ is not an addable corner ($c_1$ is the highest addable 
corner of $k+1$-residue $i$), no cell lies in the square $b$ to 
the left of $a$.  Thus, the hook of $b\wedge c_1$ is $k+1$ 
contradicting that $\gamma$ is a $k+1$-core. Therefore $x$ is not 
a $k$-hook.  
For rows corresponding to $c_j$ for $j>1$, the cells 
$x=c_{j-1}\wedge c_{j}$ have hook-length $k$ 
in $\gamma$ by Property~\ref{propodisjoint} while
cells to the right (left) of $x$ are strictly smaller (larger) 
than $k$ by Remark~\ref{hooksgrowold}.  However, since $s_i(\gamma)$ 
is obtained by adding $c_{j-1}$ and $c_j$ to $\gamma$, 
the hook of $x$ is $k+2$ in $s_i(\gamma)$ while the hooks
to the right and left of $x$ increase by one and are thus 
not $k+1$.

The proof when there is a removable corner of $k+1$-residue $i$
in $\gamma$ follows similarly.
\end{proof}

\begin{property} 
\label{sinvol}
If $\gamma\in\mathcal C_{k+1}$
then $s_i^2(\gamma)=\gamma$ for all $i\in\{0,\dots,k \}$.
\end{property} 
\begin{proof}
When there are no removable or addable corners of
$k+1$-residue $i$, $s_i$ is clearly an involution.
If $\gamma$ has at least one removable corner
of $k+1$-residue $i$ then by Property~\ref{propouter}(ii),
$s_i(\gamma)=\delta$ is the $k+1$-core where 
all removable corners $c_1,\ldots,c_n$ of $k+1$-residue $i$ have been
removed from $\gamma$.  Since Remark~\ref{both} implies 
that there can be no addable corners of $k+1$-residue $i$ in $\gamma$,
$c_1,\ldots,c_n$ are exactly the addable corners of
$\delta$ and $s_i(\delta)=\gamma$ by Property~\ref{propouter}(i).
Similar reasoning proves that $s_i$ is an involution if
$\gamma$ has at least one addable corner of $k+1$-residue $i$.
\end{proof}

In fact, the $s_i$ operators satisfy the affine Coxeter relations
(see \S~\ref{weak}).  We now conclude this section with one 
last property.

\begin{property}
\label{helpful}
For any $i=0,\ldots,k$ and $k+1$-core  $\gg$,
$s_i(\gg)$ is a $k+1$-core such that:
\begin{itemize}
\item[{(i)}] if $c_1,\dots,c_n$ is the $k+1$-string of
removable corners with $k+1$-residue $i$ in $\gamma$,
then the cells $c_1\wedge c_2,\dots,c_{n-1}\wedge c_{n}$ are the 
only cells whose hook exceeds $k$ in $\gg$ but
is $k$-bounded in $s_i(\gg)$.  
\item[{(ii)}] if $c_1,\dots,c_n$ is the $k+1$-string of
addable corners with $k+1$-residue $i$ 
in $\gamma$, then the cells
$c_1\wedge c_2,\dots,c_{n-1}\wedge c_{n}$ are the only cells 
whose hook is $k$-bounded in $\gg$ and exceeds $k$ in $s_i(\gg)$. 
\end{itemize}
\end{property}
\begin{proof}
In the case (ii) that $\gg$ has addable corners $c_1,\ldots,c_n$ of 
$k+1$-residue $i$, let $A$ denote the set of cells with hooks 
exceeding $k$ in $\gg$ and let $B$ denote such cells in 
$s_i(\gg)$.  Thus, $s\in B-A$ satisfies
$h_s(\gamma)\leq k$ and $h_s\left(s_i(\gamma) \right)>k$.
However, since $s_i(\gg)$ is a $k+1$-core by
Property~\ref{propouter} and thus has no $k+1$-hooks,
we have $h_s(\gg))\leq k$ and $h_s\bigl(s_i(\gg)\bigr)>k+1$.  
Since $s_i(\gg)$ is obtained from $\gg$ by adding corners $c_i$,
the hook of any cell $x$ in $s_i(\gg)$ has two more cells than 
the hook of $x$ in $\gg$ only if $x=c_j\wedge c_\ell$
for some $j$ and $\ell$.  However, of these, only 
$c_1 \wedge c_2,\dots,c_{n-1} \wedge c_{n}$ have a $k$-bounded hook 
in $\gg$ since the $c_i$ are separated by hooks of 
length $\ell (k+1)+k$.  This proves Case~(ii) and
Case~(i) follows by replacing $\gg \leftrightarrow s_i(\gg)$ and 
using $s_i^2=id$.
\end{proof}

\section{$k$-Young lattice and $k+1$-cores}

Recall that the set of elements covered by $\lambda$ with respect
to $\preceq$ is a subset of the partitions obtained by removing 
a corner box from $\lambda$.   These partitions must also satisfy an
additional condition that concerns the number of $k$-bounded
hooks in the $\core(\lambda)$.  Equipped with the previous
discussion of cores and their properties, 
we are now in the position to precisely understand how
the number of $k$-bounded hooks in a $k+1$-core changes
under the action of $s_i$.  This then enables us to
characterize the $k$-addable corners and consequently,
the saturated chains in the $k$-Young lattice.

\medskip

\begin{proposition}
\label{remarkone}
Given any $k$-bounded partition $\lambda$ and $\gamma=\core(\lambda)$,
\begin{equation}
\label{skeworder}
s_i\left(\gamma\right)
=
\begin{cases}
\core(\lambda- e_r)
\;
\text{where $r$ is the highest row of $\gg$ containing a removable 
corner of residue $i$}
\\
\core(\lambda+ e_r)
\;
\text{where $r$ is the highest row of $\gg$ containing 
an addable corner of residue $i$}
\\
\gg \quad\quad\quad\quad \text{when $\gamma$ has no removable 
or addable corner of residue $i$}
\end{cases}
\end{equation}
Further, when $s_i$ does not act as the identity, it acts on 
$\gg$ by removing/adding corners so that every row and column 
of $\gamma$ and $s_i(\gg)$ has the same number of $k$-bounded cells 
except in one row (and column) where $s_i(\gg)$ has one fewer/more 
$k$-bounded cell then $\gg$.  In particular, the total number of 
$k$-bounded cells in $s_i(\gamma)$ is exactly one more/fewer 
than in $\gamma$ when $\gamma$ contains an addable/removable corner 
of $k+1$-residue $i$.
\end{proposition}
\begin{proof}
Let $c_1,\dots,c_n$ be the $k+1$-string of removable corners with
$k+1$-residue $i$ in $\gamma$.
Property~\ref{helpful}(i) reveals that
the $n$ $k$-bounded cells $c_1,\dots,c_n$ in $\gg$
are not $k$-bounded in $s_i(\gg)$ while the 
$n-1$ cells $c_1 \wedge c_2,\dots,c_{n-1} \wedge c_n$
are $k$-bounded in $s_i(\gg)$ but not in $\gg$.
Therefore, $s_i$ acts on $\gg$ by decreasing the number 
of $k$-bounded cells only in the row containing
$c_1$ and in the column containing $c_n$ since 
$c_j \wedge c_{j+1}$ and $c_{j+1}$ lie in the same row 
while $c_j \wedge c_{j+1}$ and $c_{j}$ lie in the same column 
for $i=1,\dots,n-1$. Therefore, since $\lambda=\core^{-1}(\gamma)$ 
indicates the number of $k$-bounded hooks in rows of $\gamma$, we have 
$\lambda-e_r=\core^{-1}(s_i(\gamma))$ where $r$ is the highest row of $\gg$ 
containing a removable corner ($c_1$) of $k+1$-residue $i$.

Replacing $\gg \leftrightarrow s_i(\gg)$ and using $s_i^2=id$ 
proves the case with $c_1,\ldots,c_n$ the addable corners.
\end{proof}

This given, we can characterize the $k$-addable corners.

\medskip

\begin{theorem}
The order $\preceq$ can be characterized by the covering relation
\begin{equation}
\label{e1}
\lambda\cpreceq\,\mu\iff
\lambda =\mu-e_r 
\end{equation}
where $r$ is any row of $\core(\mu)$ with a removable corner 
whose $k+1$-residue $i$ does not occur in a higher removable corner, in which
case $s_i\bigl(\core(\mu) \bigr)=\core(\lambda)$.
Equivalently, $r$ can be characterized as
a row of $\core(\lambda)$ with an 
addable corner whose $k+1$-residue $i$ does not occur 
in a higher addable corner.
\label{khook}
\end{theorem}

\begin{example}  \label{whatcorners}
With $k=4$ and $\lambda=(4,2,1,1)$, 
\begin{equation}
\core(4,2,1,1) = {\tiny{\tableau*[scY]{  
\bl 1\cr 
2\cr  
3 & \bl 4\cr  
4 & 0 &  \bl 1 \cr 
 0 &  1 & 2 & 3  & 4 & 0 & \bl 1 
}}} \, ,
\end{equation}
and thus the partitions that are covered by $\lambda$ are
$(4,1,1,1)$, and $(4,2,1)$, while those that
cover it are $(4,2,1,1,1)$ and $(4,2,2,1)$.
\end{example}

\begin{proof}
Assume that $r$ is a row of $\core(\mu)$ with a removable corner $a$
of $k+1$-residue $i$.  If no removable corner of $\core(\mu)$ with 
$k+1$-residue $i$ lies higher than $a$, then Proposition~\ref{remarkone} 
implies $s_i\bigl(\core(\mu)\bigr)=\core(\mu-e_r)$, and that
the number of $k$-bounded cells of $\core(\mu-e_r)$ differs from 
$\core(\mu)$ in only one column where it is shorter by one.  Therefore, 
by the definition of $k$-conjugation, 
$(\mu-e_r)^{\omega_k} \subseteq \mu^{\omega_k}$,
implying $\mu-e_r \cpreceq\mu$. 

On the other hand, assume there is a removable corner $b$ 
of $k+1$-residue $i$ higher than $a$. To prove 
$(\mu-e_r)^{\omega_k}\not\subseteq\mu^{\omega_k}$
(implying by definition that $\mu-e_r\not\cpreceq\mu$),
it suffices to assume $b$ is a $k+1$-predecessor of $a$ since 
Property~\ref{propodisjoint} implies the removable corners 
form a $k+1$-string.  Now we shall show that a column of 
$\mu/^k$ is shorter than the same column of $(\mu-e_r)/^k$.  
The diagrams of $\mu/^k$ and $(\mu-e_r)/^k$ coincide strictly above 
row $r$ by the recursive method of constructing a $k$-skew diagram 
presented in Lemma~\ref{propo3.2} (i.e. Example~\ref{exskew}).  
In row $r$, the square $x=b\wedge a$ (see \eqref{figure}) 
must satisfy $k<h_x(\mu/^k)\leq k+1$ (or $h_x(\mu/^k)= k+1$) 
since $b$ is a $k+1$-predecessor of $a$, both removable corners. 
Therefore, deleting a cell in row $r$ allows $x\in(\mu-e_r)/^k$ without 
producing a hook exceeding $k$.  Thus the column with $x$ in $\mu-e_r$ 
is longer than the corresponding column in $\mu/^k$.
\begin{equation}
\label{figure}
{\tiny{\tableau*[scY]{  
&b \cr
& \cr
\bl & \bl x &  & & a \cr
}}} 
\quad \longrightarrow \quad  
{\tiny{\tableau*[scY]{  
&b \cr
& \cr
 & x & \bl  &\bl  & \bl a \cr
}}} 
\end{equation}

Finally, $r$ can be equivalently characterized as
the highest row in $\core(\lambda)$ with an 
addable corner of given $k+1$-residue since the 
addable corners of  $\core(\lambda)$
are exactly the removable corners of $\core(\mu)$, 
given that  $\core(\mu)=s_i\bigl( \core(\lambda)\bigr)$.
\end{proof}

Thus, we combine this result with Proposition~\ref{remarkone}
to derive the following consequences:

\begin{corollary}
\label{corokhook}
Given $k$-bounded partitions $\lambda$ and $\mu$, 
\begin{equation}
\lambda  \cpreceq\,\mu  \iff 
{\core}(\lambda)\subset {\core}(\mu)\;\;\text{and}\;\;
s_i \bigl({\core}(\lambda)\bigr) = {\core}(\mu)
\text{ for some $i \in \{0,\dots,k\}$}\,.
\end{equation}
\end{corollary}

This given, we are now able to provide a core-characterization of 
the saturated chains from the empty 
partition (hereafter $\emptyset=\lambda^{(0)}$) to any $k$-bounded 
partition $\lambda\vdash n$:
\begin{equation}
\label{defpaths}
\mathcal D^k(\lambda) = 
\left\{(\lambda^{(0)},\lambda^{(1)},\ldots,\lambda^{(n)}=\lambda)
\; : \;
\lambda^{(j)} \cpreceq\, \lambda^{(j+1)}
\right\}
\,.
\end{equation}

\begin{corollary}
The saturated chains to the vertex $\lambda\vdash n$ in the $k$-lattice 
are given by 
\label{paths}
\begin{equation*}
\mathcal D^k(\lambda)
=\left\{(\lambda^{(0)}, \lambda^{(1)}, \dots, 
\lambda^{(n)}=\lambda)
\, : \, 
{\core}(\lambda^{(j)})\subset {\core}(\lambda^{(j+1)})
\;\text{and}\;
{\core}(\lambda^{(j+1)})=s_i\left({\core}(\lambda^{(j)})\right)
\text{ for some } i
\right\}
\end{equation*}
\end{corollary}

\section{Standard $k$-tableaux}
\label{secktab}

Motivated by the proposed role of $k$-lattice chains in the study of 
certain Macdonald polynomial expansion coefficients, we pursue a tableaux
interpretation for these chains.   In this section we shall provide 
a bijection between the set of chains $\mathcal D^k(\lambda)$ and a new 
family of tableaux defined on cores.  Following our discussion in 
\S~\!\ref{weak} of the connection between the $k$-lattice and weak order 
on affine permutations, a bijection from these tableaux to certain 
reduced expressions will also be revealed.

\subsection{Definition}

\begin{definition} 
\label{defktab} 
A $k$-tableau $T$ of shape $\gg\in\mathcal C_{k+1}$ 
with $n$ $k$-bounded hooks is a filling of $\gg$ with 
integers $\{1,\ldots,n\}$ such that
 
\noindent (i) rows and columns are strictly increasing

\noindent (ii) repeated letters have the same $k+1$-residue
\end{definition}

The set of all $k$-tableaux of shape $\core(\lambda)$
is denoted by $\mathcal T^k(\lambda)$.

\begin{example}
\label{exktab}
$\mathcal T^3(3,2,1,1)$,
or the set of $3$-tableaux of shape $(6,3,1,1)$, is 
\begin{equation}
{\tiny{\tableau[scY]{7|5|4,6,7|1,2,3,4,6,7}}}
\quad \quad
{\tiny{\tableau[scY]{7|6|4,5,7|1,2,3,4,5,7}}}
\quad \quad
{\tiny{\tableau[scY]{7|4|3,6,7|1,2,4,5,6,7}}}
\quad \quad
{\tiny{\tableau[scY]{7|4|2,6,7|1,3,4,5,6,7}}}
\end{equation}
\end{example}

Our first task is to show that deleting all occurrences of the largest letter 
from a given $k$-tableau produces a new $k$-tableau.  For this we shall need 
yet another property about cores.

\begin{property} \label{propoNpluspetit}
Let $\gamma$ and $\delta$ be $k+1$-cores.  
If $\gamma\subset \delta$, then the number of $k$-bounded hooks of
$\gamma$ is smaller than that of $\delta$.
\end{property}
\begin{proof}
Let $n_\gamma$ and $n_\delta$ denote the number of 
$k$-bounded hooks of the $k+1$-cores $\gamma$ and $\delta$ respectively.
If $|\delta|=1$ then $\gamma=\emptyset$ and we have
$n_\delta=1$ while $n_\gamma=0$.  Now, assume the result holds for
all $k+1$-cores $\delta$ with $n_{\delta} < N$, and consider a pair of $k+1$-cores 
$\gamma$ and $\delta$ such that $n_{\delta}=N$ and $\gamma\subset\delta$.
For $i$ the $k+1$-residue of a removable corner of $\delta$,
$s_i(\delta)$ is a $k+1$-core whose number of $k$-bounded hooks is less than
$N$ since 
$n_{s_i(\delta)}=n_{\delta}-1$ by Proposition~\ref{remarkone}.
Thus, if $\gamma\subset s_i(\delta)$ then $n_\gamma<n_{s_i(\delta)}< N$  
by induction. 
Further, if $\gamma=s_i(\delta)$ then $n_\gamma=n_{s_i(\delta)}<N$.
Finally, in the case that $\gamma\not\subset s_i(\delta)$, 
we have $s_i(\gamma)\subseteq s_i(\delta)$ since any cell of $\gamma$
not in $s_i(\delta)$ is a removable corner of $\gamma$ of $k+1$-residue $i$.
However, $s_i^2=id$ implies that $s_i(\gamma)\subset s_i(\delta)$
and thus by induction, $n_{s_i(\gamma)}<n_{s_i(\delta)}\implies
n_{\gamma}<n_{\delta}$ by Proposition~\ref{remarkone}.
\end{proof}

\begin{definition} 
For any partition $\nu$, let $\#(\nu)$ denote 
the smallest number such that there exists a 
$k+1$-core $\gamma$ with $\#(\nu)$ $k$-bounded hooks 
and $\nu \subseteq \gamma$.
\end{definition}

\begin{lemma} \label{lemmabiz}
For any partition $\nu$, $\#(\nu)$ is the smallest number
of letters needed to fill the shape $\nu$ in such a way that

\noindent (i) rows and columns are strictly increasing

\noindent (ii) repeated letters have the same $k+1$-residue
\end{lemma}
\begin{proof}  
Let $\gamma$ be a $k+1$-core with $\#(\nu)$ $k$-bounded hooks
such that $\nu \subseteq \gamma$ and
let $n$ be the smallest number of letters needed to fill $\nu$
properly (i.e. satisfying conditions (i) and (ii)).

To show that $n\leq \#(\nu)$, it suffices to find a proper filling of
$\gamma$ using $\#(\nu)$ letters.  For $i_1$ the $k+1$-residue of a
removable corner of $\gamma$, put letter $N=\#(\nu)$ in all removable 
corners of $\gamma$ with residue $i_1$.  For $i_2$ the $k+1$-residue 
of a removable corner in $s_{i_1}(\gamma)$, put letter $N-1$
in all cells of $\gamma$ corresponding to corners of $s_{i_1}(\gamma)$ 
with residue $i_2$. By iteration, we obtain a proper filling of $\gamma$
(and consequently of its subshape $\nu$)  with $N$ letters since 
Proposition~\ref{remarkone} implies that
each $s_j$ decreases the number of $k$-bounded hooks by one.

To prove that $n\geq N$, let $T$ be a proper filling of the shape $\nu=\nu^{(n+1)}$ 
and consider the tableaux of shape $\nu^{(i)}$ obtained by
deleting the letters $i,\ldots,n$ from $T$ (for $i=n,\ldots,1$).
The lemma will follow by showing that $\#(\nu^{(i)})\geq \#(\nu^{(i+1)})-1$.
That is, starting with $\#(\nu)=N$, this would imply that
$\#(\nu^{(n)})\geq N-1$ and then $\#(\nu^{(n-1)})\geq N-2$, and
by iteration that the empty partition $\nu^{(1)}$ satisfies
$\#(\nu^{(1)})\geq N-n$.  Therefore, $0\geq N-n$.

Let $a_1,\dots,a_m$ denote the positions of the letter $n$ in $\nu$ 
and let $\bar \nu$ be the partition $\nu$ minus these cells.
It remains to show that $\#(\bar\nu)\geq \#(\nu)-1$.
Let $\bar\gamma$ and $\gamma$ denote $k+1$-cores with 
$\#(\bar\nu)$ and $\#(\nu)$ $k$-bounded hooks, respectively, and
where $\bar \nu\subseteq\bar\gamma$ and $\nu\subseteq\gamma$.
Note that since $a_1,\ldots,a_m$ are removable corners of some 
$k+1$-residue $i$ in $\nu$, they are addable corners of $\bar\nu$.
Thus, these are either addable corners of $\bar\gamma$ or
lie in $\bar\gamma$.  If all $a_1,\ldots,a_m\in\bar\gamma$
then $\nu\subseteq\bar\gamma$ implies $\#(\bar\nu)\geq \#(\nu)$
by definition of $\#(\nu)$.  Otherwise, given $a_j$ is an addable 
corner of $\bar\gamma$, the number of $k$-bounded hooks $M$ of 
$s_i(\bar\gamma)$ is $\#(\bar\nu)+1$ by Proposition~\ref{remarkone}.
However, since $\nu\subseteq s_i(\bar\gamma)$, $M\geq \#(\nu)$
and we reach our claim.
\end{proof}

\begin{proposition}
\label{si=hatsi}
Deleting all cells filled with the letter $n=|\lambda|$ 
from $T\in\mathcal T^k(\lambda)$ gives a $k$-tableau 
$\bar T\in\mathcal T^k(\mu)$, where $ {\core} (\mu)=
s_i\bigl( {\core} (\lambda)\bigr)$ for $i$ the $k+1$-residue of the cells 
containing the letter $n$.
\end{proposition}
\begin{proof}  
Let $\gamma$ be the shape of $T$, and let $\bar T$ be $T$ without letter 
$n$.  To prove that $\bar T$ is a $k$-tableau, it suffices to show that 
the shape $\nu$ of $\bar T$ is that of a $k+1$-core.  
If $i$ is the $k+1$-residue of some removable corner containing the
letter $n$ in $T$, then $s_i(\gamma)$ is a $k+1$-core with $n-1$ $k$-bounded 
hooks by Proposition~\ref{remarkone} and $s_i(\gamma)\subseteq \nu$.
If we assume by contradiction that $\nu$ is not a $k+1$-core then
$s_i(\gamma) \subset \nu$.  Thus, any $k+1$-core $\delta$ containing $\nu$ 
also satisfies $s_i(\gamma) \subset \delta$.  Therefore, 
Property~\ref{propoNpluspetit} implies that $\delta$ has more 
$k$-bounded hooks than $s_i(\gamma)$ and thus,
$\#(\nu)>n-1$ by definition.  Lemma~\ref{lemmabiz} then
leads to the contradiction saying that $\bar T$ of shape $\nu$ cannot 
be properly filled with $n-1$ letters.  Finally, since $\nu$ is
a $k+1$-core, it cannot have a removable and addable corner of the 
same $k+1$-residue by Remark~\ref{both}. Therefore, given that $i$ is 
the $k+1$-residue of a removable corner in $\gamma$ containing $n$
(thus an addable corner in $\nu$), {\it every} removable corner 
of residue $i$ in $\gamma$ contains $n$, implying $\nu=s_i(\gamma)$.
\end{proof}

\subsection{Bijection: $k$-tableaux and saturated chains}

We now introduce two maps that lead to our  bijection between chains 
$\mathcal D^k(\lambda)$ in the $k$-lattice and 
$k$-tableaux $\mathcal T^k(\lambda)$.

\begin{definition} 
For any path $P=(\lambda^{(0)},\ldots,\lambda^{(n)})
\in \mathcal D^{k}(\lambda)$, let $\Gamma(P)$ be the tableau 
constructed by putting letter $j$ in positions 
${\core}(\lambda^{(j)})\big/{\core}(\lambda^{(j-1)})$
for $j=1,\ldots,n$.

\noindent
Given $T \in \mathcal T^k(\lambda)$, let 
$\bar \Gamma(T)=(\lambda^{(0)},\ldots,\lambda^{(n)})$
where ${\core}(\lambda^{(j)})$ is the shape of the tableau
obtained by deleting letters $j+1,\ldots,n$ from $T$.
\end{definition}

To compute the action of $\Gamma$ on a path, we view the action of 
$\core$ as a composition of maps on a partition -- first skew the 
diagram and then add the squares below the skew to obtain a core.  

\begin{example}
With $k=3$: 
\begin{eqnarray*}
\left( \emptyset 
\,,
{\tiny{\tableau[scY]{|}}}
\,,
{\tiny{\tableau[scY]{,}}}
\,,
{\tiny{\tableau[scY]{,,}}}
\,,
 {\tiny{\tableau[scY]{|,,}}}
\,,
{\tiny{\tableau[scY]{||,,}}}
\,,
{\tiny{\tableau[scY]{|,|,,}}}
\,,
{\tiny{\tableau[scY]{||,|,,}}}
\;
\right)
\qquad \qquad 
&
{{\Gamma\atop\longrightarrow}\atop{\longleftarrow\atop\bar\Gamma}}
\qquad
{\tiny{\tableau[scY]{7|5|4,6,7|1,2,3,4,6,7}}}
\qquad \qquad \qquad \qquad
\qquad
\end{eqnarray*}
$$
{\rm skew} \quad\updownarrow
\qquad \qquad \qquad \qquad \qquad \qquad \qquad \qquad
$$
$$
\squaresize .2cm
\thickness .01cm
\Thickness .06cm
\emptyset \,,
\Young{\cr}
\,,
\Young{&\cr}
\,,
\Young{&&\cr}
\,,
\Young{\cr\blank&&&\cr}
\,,
\Young{\cr\cr\blank&&&\cr}
\,,
\Young{\cr&\cr\blank&\blank&&&\cr}
\,,
\Young{\cr\cr\blank&&\cr\blank&\blank&\blank&&&\cr}
\quad
{{\rm add}\atop \leftrightarrow}
\quad
\emptyset\,,
\Young{\cr }
\,,
\Young{&\cr }
\,,
\Young{&&\cr }
\,,
\Young{\cr&&&\cr }
\,,
\Young{\cr\cr&&&\cr }
\,,
\Young{\cr&\cr&&&&\cr }
\,,
\Young{\cr\cr&&\cr&&&&&\cr }
$$
\end{example}

The example suggests that $\Gamma^{-1}=\bar\Gamma$.
This will indeed follow from the following lemmas:

\begin{lemma} \label{propobi1}
If $P \in \mathcal D^k(\lambda)$, then $\Gamma(P)\in\mathcal T^k(\lambda)$.
\end{lemma}
\begin{proof}
Since the only path in $\mathcal D^k({\tiny{\tableau[sbY]{|}}}\,)$
is $P=(\emptyset,{\tiny{\tableau[sbY]{|}}})$, and 
$\Gamma(P)={\tiny{\tableau[sbY]{1}}}
\in\mathcal T^k({\tiny{\tableau[sbY]{|}}}\,)$, we proceed by
induction on $|\lambda|$.  Assume that  $\Gamma$ sends 
any path of length $n-1$ to a $k$-tableau on $n-1$ letters
and let ${  P}=(\lambda^{(0)},\ldots, \lambda^{(n)})
\in\mathcal D^k(\lambda)$.  The definition of $\Gamma$ implies that 
$\Gamma({  P})$ is obtained by adding letter $n$ to 
$T^{n-1}=\Gamma(\lambda^{(0)},\ldots, \lambda^{(n-1)})$
in positions
$\core\left(\lambda^{(n)}\right)\big/\core\left(\lambda^{(n-1)}\right)$.
Since
${\core}\left(\lambda^{(n)}\right)=s_i\left({\core}(\lambda^{(n-1)})\right)$
for some $i$ by Corollary~\ref{paths}, these positions are the addable 
corners of $\core(\lambda^{(n-1)})$ with $k+1$-residue $i$.  Therefore, 
the letters $n$ in $\Gamma(P)$ all have the same $k+1$-residue and no two
can occur in the same row or column.  Thus $\Gamma(P)$ is a
$k$-tableau of shape ${\core}(\lambda^{(n)})$ given that the
subtableau $T^{n-1}$ is a $k$-tableau by induction.
\end{proof}

\begin{lemma} \label{propobi2}
If $T \in \mathcal T^k(\lambda)$, then $\bar\Gamma(T) \in 
\mathcal D^k(\lambda)$.
\end{lemma}
\begin{proof}
Consider $T\in\mathcal T^k(\lambda)$ and let
$\bar\Gamma\left(T\right)=
(\lambda^{(0)},\ldots,\lambda^{(n)})$.
For all $j$, if $\core(\lambda^{(j)})\subset\core(\lambda^{(j+1)})$ and 
$s_i({\core}(\lambda^{(j+1)}))=\core(\lambda^{(j)})$ for some $i$,
then $\bar\Gamma\left(T\right)\in\mathcal D^k(\lambda)$ 
by Corollary~\ref{paths}.  The definition of $\bar\Gamma$ implies that
${\core}(\lambda^{(n)})$ is the shape of $T$,
and further that ${\core}(\lambda^{(n-1)})$ is the shape of
$T$ minus all occurrences of $n$.  Clearly
$\core(\lambda^{(j)})\subset\core(\lambda^{(j+1)})$,
and further
$s_i\bigl({\core}(\lambda^{(n)})\bigr)={\core}(\lambda^{(n-1)})$
where $i$ is the $k+1$-residue of the cells containing $n$
by Proposition~\ref{si=hatsi}.  Thus, the codomain of $\bar\Gamma$ 
is $\mathcal D^k(\lambda)$ by iteration.
\end{proof}

We are now set to prove that $\Gamma$ is a bijection between 
saturated chains and $k$-tableaux.

\begin{theorem} \label{firstbijection}
$\Gamma$ is a bijection between 
$\mathcal D^k(\lambda)$ and $\mathcal T^k(\lambda)$ with
$\Gamma^{-1}=\bar\Gamma$. 
\end{theorem}
\begin{proof}
From Lemmas~\ref{propobi1} and \ref{propobi2}, it suffices to
prove that $\Gamma$ and $\bar \Gamma$ are inverses.
We start by showing that $\bar\Gamma\Gamma({P})={P}$.
Given $\Gamma(\lambda^{(0)},\ldots,\lambda^{(n)})=T$,
we must show that if
$\bar\Gamma\left(T\right)=(\mu^{(0)},\ldots,\mu^{(n)})$
then ${\core}(\mu^{(\ell)})={\core}(\lambda^{(\ell)})$ 
for $\ell=0,\ldots,n$.  The definition of $\Gamma$ implies that
shape$(T)={\core}(\lambda^{(n)})$ and that the letter $n$ lies 
in ${\core}(\lambda^{(n)})\big/{\core}(\lambda^{(n-1)})$.
At the same time, the definition of $\bar\Gamma$
implies that shape$(T)={\core}(\mu^{(n)})$
and ${\core}(\mu^{(n-1)})$ is the shape of the tableau
obtained by deleting all occurrences of the letter $n$ from $T$.
Therefore, ${\core}(\mu^{(n-1)})={\core}(\lambda^{(n-1)})$.
By iteration, $\bar\Gamma\Gamma({  P})={  P}$.

On the other hand, given $\bar\Gamma(T)=(\lambda^{(0)},\ldots,\lambda^{(n)})$,
we must show that if $\Gamma(\lambda^{(0)},\ldots,\lambda^{(n)})=\bar T$
then $\bar T=T$.  The definition of $\Gamma$ implies
that $\bar T$ is the tableau obtained by filling the cells of
${\core}(\lambda^{(j+1)})/{\core}(\lambda^{(j)})$
with letter $j+1$.  However, by definition of $\bar\Gamma$, 
${\core}(\lambda^{(j)})$ is the shape 
obtained by deleting the letters $j+1,\dots,N$ from $T$, and
${\core (\lambda^{(j+1)})}$ is the shape obtained by deleting 
$j+2,\dots,N$ from $T$.   Therefore, the cells
${\core}(\lambda^{(j+1)})/{\core}(\lambda^{(j)})$ in $T$
are filled with letter $j+1$.
Thus $T=\bar T$.
\end{proof}

\section{The $k$-Young lattice and the 
weak order on $\tilde S_{k+1}/S_{k+1}$} \label{weak}

In this section we shall see how the $k+1$-core characterization
of the $k$-Young lattice covering relations given in
Corollary~\ref{corokhook} leads to the identification of the 
$k$-Young lattice as the weak order on $\tilde S_{k+1}/S_{k+1}$.  
A by-product of this result is a simple bijection between
reduced words and $k$-tableaux and one between
$k$-bounded partitions and affine permutations in $\tilde S_{k+1}/S_{k+1}$.

\subsection{The isomorphism}
To establish that the $k$-Young lattice is isomorphic to 
weak order on the set of minimal coset representatives of 
$\tilde S_{k+1}/S_{k+1}$, we rely foremost on the fact 
\cite{[L]} that the $s_i$ operators satisfy the 
affine Coxeter relations \eqref{coxeter}, and thus provide a realization 
of the affine symmetric group on $k+1$-cores.
\begin{property}
 The $s_i$ operators satisfy 
\begin{eqnarray} \label{coxetersi}
 s_i^2  =  id, \qquad
 s_i  s_j =  s_j  s_i \quad \; (i-j\neq \pm 1\mod k+1),\quad\text{and}\quad
 s_i  s_{i+1}  s_{i} =   s_{i+1}  s_i  s_{i+1} \, .
\end{eqnarray}
\end{property}

\noindent The following map is then well defined:

\begin{definition} 
For $\sigma\in\tilde S_{k+1}$, let $\shpe$ send
$\sigma$ to a $k+1$-core by
\begin{equation}
\shpe: \sigma =s_{i_1} \cdots s_{i_\ell}\cdot\emptyset\,,
\end{equation}
where $i_1\cdots i_\ell$ is any reduced word for $\sigma$
and $\emptyset$ is the empty $k+1$-core.
\end{definition}

A characterization for Bruhat order in terms of the containment
of cores stemming from this map is provided by Lascoux 
in \cite{[L]}.  To be precise,

\begin{proposition}
\label{proplascoux} 
The map $\shpe: \tilde S_{k+1}/S_{k+1} \to {\mathcal C}_{k+1}$ is 
an isomorphism from Bruhat order on $\tilde S_{k+1}/S_{k+1}$ to Young 
order ($\subseteq$) on ${\mathcal C}_{k+1}$.
\end{proposition}

We are thus able to obtain from our $k+1$-core characterization of 
the chains in the $k$-lattice that this lattice is isomorphic 
to the weak order on $\tilde S_{k+1}/S_{k+1}$:

\begin{corollary}
\label{corowords}
Let $\sigma,\tau\in\tilde S_{k+1}/S_{k+1}$,
and let $\lambda=\kbnd(\shpe(\sigma))$ and $\mu=\kbnd(\shpe(\tau))$.
Then
\begin{eqnarray}
\sigma<\!\!\!\!\cdot_{w}\,\tau \quad \iff \quad \lambda\cpreceq\mu \, .
\end{eqnarray}
\end{corollary}
\begin{proof}
Proposition~\ref{proplascoux} implies a characterization 
of the covering relations for weak order on 
$\tilde S_{k+1}/S_{k+1}$.  That is, since $\shpe$ is a bijection
and the weak order is a suborder of the Bruhat order, we have
for $\sigma,\tau\in\tilde S_{k+1}/S_{k+1}$
\begin{equation}
\label{wklascoux}
\sigma<\!\!\!\!\cdot_{w}\,\tau
\quad \iff \quad \shpe(\sigma)\subset \shpe(\tau) \quad \text{\rm ~and } \quad
s_i \, \shpe(\sigma)=\shpe(\tau) \text{\rm ~for some } i
\,.
\end{equation}
The result thus follows from the characterization of
$\rightarrow_k$ given in Corollary~\ref{corokhook}.
\end{proof}

\subsection{Bijection: $k$-tableaux and reduced words} 
We have seen in Theorem~\ref{firstbijection} that the
saturated chains to shape $\lambda$ in the $k$-lattice 
are in bijection with $k$-tableaux of shape $\kbnd(\gamma)$.
On the other hand, the reduced words for $\sigma\in\tilde S_{k+1}/S_{k+1}$ 
encode the chains to $\sigma$.  Corollary~\ref{corowords} thus
implies there is a bijection between $k$-tableaux
of shape $\gamma$ and the reduced words for $\shpe^{-1}(\gg)$.

This bijection arises naturally by noting from Corollary~\ref{paths} 
that the association between a $k$-tableau and chain 
$(\lambda^{(0)},\dots,\lambda^{(n)}=\lambda)$ in the $k$-Young lattice 
is determined by a sequence $i_n\cdots i_2\,i_1$ such that 
$s_{i_j}(\core(\lambda^{(j-1)})) =\core(\lambda^{(j)})$ for $j=1,\dots,n$. 
However, this sequence can also be viewed as a reduced word 
for the permutation $\sigma$ where $\shpe(\sigma)=\core(\lambda)$
by Eq.~\eqref{wklascoux}.  Therefore, the following map provides the 
desired bijection:

\begin{definition}
For a $k$-tableau $T$ with $m$ letters
where $i_{a}$ is the $k+1$-residue of the letter $a$, define
\begin{equation*}
{\mathfrak w}: T\mapsto {i_m} \cdots {i_1}
\,.
\end{equation*}

\noindent
For $w={i_m}\cdots {i_1} \in Red(\sigma)$, 
$\mathfrak w^{-1}(w)$ is the tableau with letter $\ell=1,\ldots,m$ 
occupying the cells of $s_{i_\ell} \cdots s_{i_1} \cdot \emptyset\big/
s_{i_{\ell-1}} \cdots s_{i_1} \cdot \emptyset$.
\end{definition}

\begin{example}
With $k=3$:
$$
T\,=\,
{\tiny{\tableau[scY]{7|6|4,5,7|1,2,3,4,5,7}}}
\qquad
{\mathfrak w\atop \leftrightarrow}
\qquad
\, 1\,2\,0\,3\,2\,1\,0
\quad
\text{ since the $4$-residues are }
\quad
{\tiny{\tableau[scY]{1|2|3,0,1|0,1,2,3,0,1}}}
$$
\end{example}

\begin{proposition} \label{theoword}
\label{bijecTabword}
The map ${\mathfrak w} : \mathcal T^k({\lambda})  
\longrightarrow Red(\sigma)$ is a bijection,
where $\sigma\in \tilde S_{k+1}/S_{k+1}$ is defined
uniquely by $\core(\lambda)=\shpe(\sigma)$.
\end{proposition}

We will now make use of canonical chains in the $k$-Young lattice
to obtain a simple bijection between $k$-bounded partitions and
permutations in $\tilde S_{k+1}/S_{k+1}$.

\begin{definition}
For any partition $\lambda$, let ``$w_{\lambda}"$
be the word obtained by reading
the $k+1$-residues in each row of $\lambda$, from right to left,
starting with the highest removable corner and ending in the first cell
of the first row.  Further, let ``$\sigma_{\lambda}$" be the affine 
permutation corresponding to $w_{\lambda}$. 
\end{definition}

\begin{example} \label{exampleres}
For $\lambda=(3,2,2,1)$ and $k=3$, $w_\lambda=1\, 3\,2\,0\,3\,2\,1\,0$ 
and $\sigma_{\lambda}=
\hat s_{1}\hat s_{3}\hat s_{2}\hat s_{0}
\hat s_{3}\hat s_{2}\hat s_{1} \hat s_{0}$
since:
\begin{equation}
\lambda={\tiny{\tableau[scY]{1|2,3|3,0|0,1,2}}}
\end{equation}
\end{example}

\begin{proposition}  $\sigma_{\lambda}$ belongs to 
$\tilde S_{k+1}/S_{k+1}$ and 
$\shpe(\sigma_{\lambda})=\core(\lambda)$.
\end{proposition}
\begin{proof}
Consider $\lambda\in\mathcal P_k$.  In light of Proposition~\ref{theoword}, 
it suffices to show that there is some $k$-tableau $T$ of shape 
$\core(\lambda)$ where ${\mathfrak w}(T)=w_\lambda$.  Note that
Corollaries~\ref{paths} and \ref{firstbijection} imply $\mathfrak w(T)$ 
(of shape $\core(\lambda)$) is obtained from a certain
chain $(\lambda^{(0)},\dots,\lambda^{(n)}=\lambda)$ in the $k$-Young 
lattice by taking the sequence $i_n\cdots i_2\,i_1$ such that 
$s_{i_j}(\core(\lambda^{(j-1)})) =\core(\lambda^{(j)})$ for $j=1,\dots,n$. 
Now, there exists a canonical saturated chain $P_\lambda$ (and thus a 
canonical sequence $i_n\cdots i_2\,i_1$) such that
$\lambda^{(j)}$ is obtained by removing the highest removable corner 
of $\lambda^{(j+1)}$.  The existence of such a chain is ensured by
Theorem~\ref{khook} since the highest removable corner of a 
partition is always the highest of its $k+1$-residue.  However,
the highest removable corner of a partition $\lambda$ coincides with 
the highest removable corner of $\core(\lambda)$ and we therefore find
that $i_n\cdots i_2\,i_1$ is exactly $w_\lambda$.
\end{proof}

Given the bijection between $k$-bounded partitions and $k+1$-cores,
this immediately provides a bijection between $k$-bounded partitions and 
permutations in $\tilde S_{k+1}/S_{k+1}$.

\begin{corollary}
The map $\phi:\mathcal P_k\to\tilde S_{k+1}/S_{k+1}$ where 
$\phi(\lambda)=\sigma_\lambda$ is a bijection whose inverse
is $\phi^{-1}=\kbnd\circ\shpe$.
\end{corollary}

\begin{example}
Given $\sigma\in\tilde S_4^I$ with a reduced expression
$w=3\,1\,0\,3\,2\,1\,3\,0$,
we construct the shape:
\begin{equation*}
s_3s_1 s_0 s_3 s_2 s_1s_3s_0
\cdot \emptyset =
{\tiny{\tableau[scY]{0}}}
\;\to\;
{\tiny{\tableau[scY]{3|0}}}
\;\to\;
{\tiny{\tableau[scY]{3|0,1}}}
\;\to\;
{\tiny{\tableau[scY]{2|3|0,1,2}}}
\;\to\;
{\tiny{\tableau[scY]{2|3|0,1,2,3}}}
\;\to\;
{\tiny{\tableau[scY]{2|3,0|0,1,2,3,0}}}
\;\to\;
{\tiny{\tableau[scY]{1|2|3,0,1|0,1,2,3,0,1}}}
\;\to\;
{\tiny{\tableau[scY]{1|2,3|3,0,1|0,1,2,3,0,1}}}
\end{equation*}
from which we read the number of $3$-bounded hooks to obtain
$\phi^{-1}(\sigma)=(3,2,2,1)$.  Conversely, $\sigma$ can be recovered
from $(3,2,2,1)$ by using Example~\ref{exampleres} to find $\phi(3,2,2,1)
=\hat s_{1}\hat s_{3}\hat s_{2}\hat s_{0}\hat s_{3}\hat s_{2}\hat s_{1} \hat s_{0}$
(one easily checks that 
$ 3\,1\,0\,3\,2\,1\,3\,0$ and $1\,3\,2\,0\,3\,2\,1\,0$ are reduced
words for the same permutation). 
\end{example}

\medskip

\section{Comparing elements differing by more than one box}

Now that we have been able in \S~\!\ref{secktab}
to explicitly understand the covering relation 
for the $k$-order and to characterize the chains, it
is natural to ask what can be said about the relation among vertices 
differing by more than one box.  In this section we shall prove that
\begin{quote}
If $\mu/\lambda$ and $\mu^{\omega_k}/\lambda^{\omega_k}$ are
horizontal and vertical strips respectively, then
$\lambda \preceq \mu$. 
\end{quote}
A number of somewhat technical properties will lead us
to this result and shall also be used in our development of
a semi-standard version of the $k$-tableaux corresponding to 
certain chains in the $k$-Young lattice.
We begin by continuing the study of $k+1$-cores, 
concentrating on pairs $\gamma\subseteq\delta$.

\begin{definition}  
Let $\gamma$ and $\delta$ be $k+1$-cores with
$\gamma\subseteq\delta$.  A {\it ``rowadder"} is a 
cell $s\in\delta/\gamma$ such that there is no 
cell in $\dg$ that is a $k+1$-predecessor of $s$.
\end{definition}  

Two properties concerning the existence of rowadders 
are needed.

\begin{property} \label{twoinone}
If $\gamma$ and $\delta$ are $k+1$-cores with $\gamma\subseteq\delta$, 
then $\delta/\gamma$ has a rowadder at the top of the leftmost 
column that contains more than one cell.
\end{property}
\begin{proof}
Let $b$ (of $k+1$-residue $i$) denote the cell in $\dg$ at the top 
of the leftmost column with more than one cell.  
Note that if $x\in\gamma$ lies immediately southwest of $b$,
then no cell of $\gamma$ lies to the right of $x$.
Further, the diagram of $\delta/\gamma$ to the left of $b$
is a series of horizontal rows since columns to the left of $b$ 
have at most one cell. 
\begin{equation}
{\tiny{\tableau*[scY]
{&\bar b&&\cr\bl&\bl&\bl&\bl & &&& &\bl &\bl \cr
\bl &\bl &\bl &\bl  &\bl&\bl &\bl &\bl &\bl  &\bl&&&b&
\cr
\bl &\bl &\bl &\bl &\bl&\bl &\bl &\bl &\bl & \bl & \bl & \bl x&&&& }}} 
\end{equation}
Suppose by contradiction that there is a cell $\bar b\in\delta/\gamma$ 
that is a $k+1$-predecessor of $b$.  Then $h_{\bar b\wedge x}(\gamma)=k+1$,
violating the assumption that $\gamma$ is a $k+1$-core.
\end{proof}

\begin{remark} \label{horistrip}
For partitions $\lambda$ and $\mu$,
$\mu/\lambda$ is a horizontal strip iff $\lambda \subseteq \mu$
and $\lambda_r \geq \mu_{r+1}$ for all $r$.  Further,
$\mu/\lambda$ is a vertical strip iff $\mu_r-\lambda_r \in \{0,1\}$
for all $r$.
\end{remark}

\medskip

\begin{property} \label{combine} 
Consider $\gamma=\core(\lambda)$ and 
$\delta=\core(\mu)$ with $\gamma\subseteq\delta$.
Let $\ell$ denote the leftmost column of $\dg$ with
more than one cell.

\noindent
(i) If there are rowadders
in the top two cells of $\dg$ in column $\ell$,
then $\mu/\lambda$ is not a horizontal strip. 

\noindent
(ii) If there is no rowadder in the second row of $\dg$ in
of column $\ell$,
then $\lambda^{\omega_k} \not \subseteq \mu^{\omega_k}$.  
\end{property} 
\begin{proof}
Case (i):  the number of $k$-bounded hooks in row  $r$ of $\delta$ 
(resp. $\gamma$) is $\mu_r$ (resp. $\lambda_r$).  
Thus, by Remark~\ref{horistrip}, it suffices to prove that 
there are at least $\lambda_r+1$ $k$-bounded hooks in row $r+1$ 
of $\delta$ for some $r$.
We shall consider the rows $r+1$ and $r$ containing rowadders 
$a,b\in\delta/\gamma$.  If $\bar a$ denotes the extremal cell 
of $\gamma$ that is a $k+1$-predecessor of $a$, then 
the extremal cell $\bar b$ of $\gamma$ that $k+1$-precedes
$b$ either lies below or beside $\bar a$ 
since $\bar a$ is extremal.  However, if $\bar b$ lies beside $\bar a$,
the hook of $\bar b\wedge b$ is $k+1$ in the $k+1$-core $\gamma$
implying this case does not occur.
When $\bar b$ lies below $\bar a$, the square $\hat b$ to the right of $\bar a$ is not 
in $\delta$ since $\bar b$ is extremal in $\gamma$ and $b$ is a rowadder:  
\begin{eqnarray}
{\scriptstyle{\tableau[scY]{
\tf,,,
|\tf,\tf,\tf,\tf\bar a,\bl \hat b,\bl,\bl,\bl \hat c
|\tf,\tf,\tf,\tf\bar b,\tf,\tf,\tf
|\tf,\tf,\tf,\tf,\tf,\tf,\tf,\tf,\tf
|\tf,\tf,\tf,\tf x_a, \tf,\tf ,\tf ,\tf x_c,\tf,,,a,,,,c ,\bl,\bl,
\bl\leftarrow,\bl, 
\bl r+1
|\tf,\tf,\tf,\tf y_b,\tf x_b,\tf,\tf,\tf,\tf,\tf,\tf,b,,,,, ,\bl,
\bl\leftarrow,\bl, \bl  r
}}}
\end{eqnarray}
Notice that the hook of $x_b=\hat b\wedge b$ in $\gamma$ is
$k$-bounded while the hook of $y_b=\bar b\wedge b$ in $\gamma$ 
exceeds $k$.  Therefore,  $\lambda_r$ is the number of cells 
strictly between $y_b$ and $b$ (equivalently, $x_a$ and $a$).
To determine the number of $k$-bounded hooks of $\delta$,
let $c$ denote the last cell in row $r+1$ of $\delta$ 
and $\hat c$ the square a $k+1$-predecessor of $c$ in the row
with $\hat b$.  Since $\hat c$ does not belong to $\delta$,
the hook length of $x_c=\hat c \wedge c$ is at most $k+1$.  But because
$\delta$ is a $k+1$-core, $x_c=\hat c \wedge c$ thus has a $k$-bounded hook 
in $\delta$ as do all the cells of $\delta$ to the right of $x_c$.
Given that the number of cells strictly between $x_c$ and $c$
equals the number of cells, $\lambda_r$, strictly between $x_a$ and $a$, 
we have at least $\lambda_r+1$ $k$-bounded hooks in row $r+1$ of $\delta$
as claimed.

Case (ii): Since the number of $k$-bounded hooks in a column of 
$\gamma$ (resp. $\delta$) corresponds to a row of 
$\lambda^{\omega_k}$ (resp. $\mu^{\omega_k}$),
to prove $\lambda^{\omega_k} \not \subseteq \mu^{\omega_k}$,
it suffices to show that there are more $k$-bounded 
hooks in some column of $\gamma$ than in that column
of $\delta$.  Let $a$ (of $k+1$-residue $i$) denote the top cell 
in the first column $\ell_a$ of 
$\delta/\gamma$ containing more than one cell.  
By assumption, the cell $b$ below $a$ is not a rowadder
and thus there is a cell $\bar b\in\delta/\gamma$ of the same 
$k+1$-residue as $b$ to the left of column $\ell_a$. 
Hence the square $\bar a$ above $\bar b$ has $k+1$-residue $i$
and since the diagram of $\delta/\gamma$ to the left of column $\ell_a$
is a series of horizontal rows,
$\bar a\not\in\delta/\gamma$. 
\begin{equation}
{\tiny{\tableau*[scY]
{
\bl&\bl \bar a\cr
&\bar b&&\cr\bl&\bl&\bl&\bl & &&& & &\bl \cr
\bl &\bl x_a &\bl &\bl &\bl&\bl &\bl &\bl &\bl  &&&&a&
\cr
\bl &\bl x_b &\bl &\bl &\bl&\bl &\bl &\bl &\bl & \bl & \bl & \bl&b&&& }}} 
\end{equation}
Since $\bar a$ and $a$ have the same $k+1$-residue,
the cell $x_a=\bar a\wedge a$ has hook-length bounded by $k$ in $\gamma$ 
and at least $k+1$ in $\delta$.  Similarly for the cell $x_b=\bar b\wedge b$.
Therefore, in the column with $x_a$, $\bar b\in \delta$ is the only 
$k$-bounded hook in $\delta$ that is not in $\gamma$ while $x_a$ 
and $x_b$ are $k$-bounded hooks 
in $\gamma$ that are not $k$-bounded in $\delta$. We reach our claim since
$\gamma$ has at least one more $k$-bounded hook 
than $\delta$ in this column.
\end{proof}

We shall say that $\mu,\lambda$ are {\it admissible} iff $\mu/\lambda$ and
$\mu^{\omega_k}/\lambda^{\omega_k}$ are respectively horizontal 
and vertical strips, {\it i.e.} iff $\mu,\lambda$ are $r$-admissible 
for some $r$.

\begin{proposition} \label{propo64}
If $\mu,\lambda\in\mathcal P_k$ forms an admissible pair,
then ${\core}(\mu)/{\core}(\lambda)$ is a horizontal 
strip. 
\end{proposition}
\begin{proof} 
Given $\mu,\lambda$ is an admissible pair, we have
$\mu/\lambda$ is a horizontal strip and $\mu^{\omega_k}/\lambda^{\omega_k}$ 
is a vertical strip.  In particular, $\lambda \subseteq \mu$ and
by Property~\ref{onedirect}, $\core(\lambda)\subseteq \core(\mu)$.  
Now, assume by contradiction that $\core(\mu)/\core(\lambda)$ contains 
some column with more than one cell.  The top cell $c$ of the leftmost 
such column must be a rowadder by Property~\ref{twoinone}.  If the 
cell $\bar c$ below $c$ is a rowadder, then this column contains two 
rowadders implying by Property~\ref{combine}(i) that $\mu/\lambda$ is not a 
horizontal strip.  On the other hand, if $\bar c$ is not
a rowadder, then $\lambda^{\omega_k}\not\subseteq\mu^{\omega_k}$ 
by Property~\ref{combine}(ii) and thus $\mu^{\omega_k}/\lambda^{\omega_k}$  
is not a vertical strip.  Either case gives a contradiction.
\end{proof}

\medskip

\begin{lemma}  \label{lemrtmost}
Let $\gamma$ and $\delta$ be $k+1$-cores where no column 
has more $k$-bounded hooks in $\gg$ than in $\delta$,
and where $\delta/\gamma$ is a horizontal strip.
With $i$ denoting the $k+1$-residue of the rightmost cell in $\delta/\gamma$,
the removable corners of $k+1$-residue $i$ in $\delta$ 
are exactly the cells of $k+1$-residue $i$ in $\delta/\gamma$.
\end{lemma}
\begin{proof} 
Let $a_1$ (of $k+1$-residue $i$) denote the rightmost cell 
in $\delta/\gamma$ and note that $a_1$ is a removable corner 
since $\delta/\gamma$ is a horizontal strip.  If $a_1$ is not 
a rowadder of $\dg$, then there is a cell $a_2\in\delta/\gamma$ 
that is a $k+1$-predecessor of $a_1$.  Similarly, if $a_2$ is 
not a rowadder then there is a cell $a_3\in\dg$ which is a $k+1$-predecessor 
of $a_2$.  By iteration, we eventually reach a rowadder $a_m\in\dg$,
and have the $k+1$-string $a_1,a_2,\ldots,a_m$ of cells with 
$k+1$-residue $i$.  Note that $a_1,\ldots,a_m$ are all extremal cells
of $\delta$ since they lie in the horizontal strip $\dg$.  Furthermore, 
no cell lies to the right of $a_1$ implying that no cell lies to the 
right of any extremal cell with $k+1$-residue $i$ above $a_1$, by 
Property~\ref{propo4.1}.  Therefore, $a_1,\ldots,a_m$ are all removable 
corners of $\delta$.  It thus remains to show that any extremal cell of 
$k+1$-residue $i$ in $\delta$ above $a_m$ or below $a_1$ is not removable.

The diagrams of $\gamma$ and $\delta$ coincide south-east of $a_1$,
given $a_1$ is the rightmost element of $\dg$.  If $a_1$ is a
$k+1$-predecessor of an extremal cell $d$, then a cell must
lie to the right of $d$ since otherwise, $h_{a_1\wedge d}(\gg)=k+1$ 
in the $k+1$-core $\gamma$.  Property~\ref{propo4.1} 
thus implies that all extremal cells of $k+1$-residue $i$ 
lying south-east of $d$ also have a cell to their right.  
Therefore there are no removable corners of $k+1$-residue $i$  
south-east of $a_1$.

Similarly by Property~\ref{propo4.1}, our claim will follow
by showing that there is a cell $b\in\delta$ above the 
extremal cell $a_{m+1}\in\delta$ that is a $k+1$-predecessor of 
$a_m$.
Suppose $b\not\in\delta$.  Then the hook length of 
$a_{m+1}\wedge a_m$ is $k+2$ in $\delta$ since $a_{m+1}$ and 
the removable corner $a_m$ have the same $k+1$-residue, but is
$k$-bounded in $\gamma$ since $a_m\not\in\gamma$ 
and $\gamma$ has no $k+1$-hooks.  
Note that the column containing $a_{m+1}$ is of the same length in $\gamma$ 
as in $\delta$ since $b\not\in\delta$ and $a_{m+1}\not\in\dg$.
Therefore, $\gamma$ has more $k$-bounded hooks in this 
column  contradicting our assumption.
\end{proof} 

\begin{theorem} \label{coro67}
If $\mu,\lambda$ is $n$-admissible, then there are
distinct integers $i_1,\dots,i_n$ where
$$
{\core}(\lambda)=s_{i_1}\cdots s_{i_n}\bigl({\core}(\mu)\bigr)\,.
$$
\end{theorem}
\begin{proof}
Since $\mu^{\omega_k}_r$ is the number of $k$-bounded hooks in 
column $r$ of $\delta={\core}(\mu)$, given that $\mu,\lambda$ is 
$n$-admissible, no column has more $k$-bounded hooks in 
$\gamma={\core}(\lambda)$ than in $\delta={\core}(\mu)$.
Further, $\dg$ is a horizontal strip by Proposition~\ref{propo64}.  
Therefore, if $i_N$ denotes the $k+1$-residue of the rightmost cell 
$a_N\in\dg$, then Lemma~\ref{lemrtmost} implies that the diagram 
$s_{i_N}(\delta)/\gamma$ can be obtained by deleting all cells of 
$k+1$-residue $i_N$ from $\dg$ and is thus a skew
diagram with no more than one cell in each column.

We now claim that no column has more $k$-bounded hooks in $\gamma$ than in 
$s_{i_N}(\delta)$.  Proposition~\ref{remarkone} gives 
that $s_{i_N}(\delta)$ has 
the same number of $k$-bounded hooks as $\delta$ in every column except
the one containing the cell $a_N$,  where it has one fewer.  
Since no column has more $k$-bounded hooks in $\gamma$ than in $\delta$,
it suffices to show that in the column with $a_N$, $\gamma$ does not have 
more $k$-bounded hooks than $s_{i_N}(\delta)$.
This follows by noting that weakly to the right of the
column with $a_N$, $s_{i_N}(\delta)$ and $\gamma$ coincide.

Therefore we can use Lemma~\ref{lemrtmost} to prove that
$s_{i_{N-1}}(s_{i_N}(\delta))/\gamma$ can be obtained by 
deleting all cells of $k+1$-residue $i_N$ and $i_{N-1}$ from 
$\delta/\gamma$.  By iterating the preceding argument,
there is some $N$ where $s_{i_1}\cdots s_{i_N}(\delta)/\gamma$ 
is the empty partition implying that $\gamma=s_{i_1}\cdots s_{i_N}(\delta)$.
Since each iteration causes the removal of all cells with a given 
$k+1$-residue from $\dg$, $i_1,\ldots,i_N$ are distinct. Further,
since the number of $k$-bounded hooks in $\delta$ is lowered by one 
each time by Proposition~\ref{remarkone}, $N=|\mu|-|\lambda|=n$.
\end{proof}

Using this result, Corollary~\ref{corokhook} implies

\begin{corollary} \label{coronouveau}
If $\mu,\lambda$ is an admissible pair,
then $\lambda \preceq \mu$. 
\end{corollary}

We conclude this section with another set of conditions 
under which $\lambda\preceq\mu$.  

\bigskip

\begin{theorem}  \label{lemhoriverti}
If $\lambda\subseteq\mu$, $\lambda^{\omega_k}\subseteq\mu^{\omega_k}$, 
and ${\core}(\mu)/{\core}(\lambda)$ is a horizontal strip,
then $\mu,\lambda$ is admissible.
\end{theorem}

\begin{proof}
We start by showing that $\mu/\la$ is a horizontal strip, or
equivalently by Remark~\ref{horistrip}, that the number of $k$-bounded 
hooks in row $r$ of $\gamma={\core}(\lambda)$ is not smaller than the 
number of $k$-bounded hooks in row $r+1$ of $\delta={\core}(\mu)$.

In row $r$ of $\gamma$, let $y_r$ denote the last cell and let $x_r$ 
be the rightmost cell with a hook exceeding $k$.  Note that 
$h_{x_r}(\gamma)>k+1$ since $\gamma$ is a $k+1$-core.  If there
are $d-1$ cells strictly between $x_r$ and $y_r$, then $\gamma$
has $d$ $k$-bounded hooks in row $r$.  It thus remains to prove that 
there are no more than $d$ $k$-bounded hooks in row $r+1$ of $\delta$.
In row $r+1$ of $\delta$, let $y_a$ be the last cell and let $x_a$ 
be the cell so that there are $d-1$ cells between $x_a$ and $y_a$ 
(if $x_a\not\in\delta$ then $\delta_{r+1}\leq d$ and the claim holds).  
$$
\squaresize .4cm
\thickness .01cm
\Thickness .06cm
\Young{ udrl&&
\cr
udrl&udrl&udrl&udrl&
\cr udrl &udrl&udrl&udrl&udrl&udrl&udrl&udrl&
\cr udrl&udrl&udrl&udrl {\scriptstyle x_a}& udrl&udrl&udrl&udrl&udrl&&& 
{\scriptstyle y_a}
\cr udrl&udrl&udrl&udrl&udrl 
{\scriptstyle x_r}&udrl&udrl&udrl&udrl&udrl&udrl&udrl&
udrl {\scriptstyle y_r}&&& &&\blank $\leftarrow$ {\rm row } $r$\cr }$$
Note that $x_a$ lies weakly to the left of $x_r$ since $y_a$ 
lies weakly to the left of $y_r$ given $\dg$ is a horizontal strip.  
Thus, the number of cells above $x_a$ in $\gamma$ 
is weakly greater than $\ell_r-1$ for $\ell_r$ the number of cells above 
$x_r$ in the partition $\gamma$.
Since $\gamma\subseteq\delta$, the number of cells $\ell_a$
above $x_a$ in $\delta$ satisfies $\ell_a \geq \ell_r-1$.
Hence, $h_{x_a}(\delta)=\ell_a+d+1 \geq \ell_r+d=h_{x_r}(\gamma)-1>k$.
That is, $h_{x_a}(\delta)$ exceeds $k$.
Therefore the maximal number of $k$-bounded hooks in row $r+1$ of 
$\delta$ is $d$.

To see that $\mu^{\omega_k}/\lambda^{\omega_k}$ is a vertical strip,
note that $\delta/\gamma$ has at most
one box in every column.  Thus, the number of $k$-bounded hooks in
a column of $\delta$ cannot exceed the number of
$k$-bounded hooks in that column of $\gamma$ by more than one
since $\gamma \subseteq \delta$ implies any hook exceeding $k$ in
$\gamma$ must exceed $k$ in $\delta$.  Now, recall that the number of 
$k$-bounded hooks in the columns of $\gamma$ and $\delta$ are 
$\lambda^{\omega_k}$ and $\mu^{\omega_k}$ respectively.  
Given $\lambda^{\omega_k} \subseteq \mu^{\omega_k}$, this leads to
$\mu_r^{\omega_k} -\lambda_r^{\omega_k} \in \{ 0,1\}$ for all $r$
-- conditions that are equivalent to $\mu^{\omega_k}\big/ \lambda^{\omega_k}$
being a vertical strip.
\end{proof}

\section{Generalized $k$-tableaux and the $k$-Young lattice}

We now introduce a set of tableaux that serve 
as a semi-standard version of $k$-tableaux.

\begin{definition}
\label{defktabgen} 
Let $\gg$ be a $k+1$-core, $m$ be the number of $k$-bounded hooks 
of $\gg$, and $\aa=(\aa_1,\ldots,\aa_r)$ be a composition of
$m$. A semi-standard $k$-tableau of shape $\gg$ and evaluation $\aa$ 
is a filling of $\gg$ with integers $1,2,\ldots,r$ such that

\noindent
(i) rows are weakly increasing and columns are strictly increasing

\noindent
(ii) the collection of cells filled with letter $i$ are labeled with exactly 
$\alpha_i$ distinct $k+1$-residues.
\end{definition}

We denote the set of all semi-standard $k$-tableaux of shape 
$\core(\lambda)$ and evaluation $\aa$ by $\mathcal T^k_{\aa}(\lambda)$.  
When $\aa=(1^m)$, we call the $k$-tableaux ``{\it standard}''.
In this case, $\mathcal T^k_{(1^m)}(\lambda)$ is the set 
$\mathcal T^k(\lambda)$ of $k$-tableaux introduced in \S~\!\ref{secktab}.  
Hereafter, a semi-standard $k$-tableau will simply be
referred to as a $k$-tableau.

\begin{example}
\label{exssktab}
For $k=3$, $\mathcal T_{(1,3,1,2,1,1)}^{3}(3,3,2,1)$ of 
shape $\core\left((3,3,2,1)\right)=(8,5,2,1)$ is the set:
\begin{equation}
{\tiny{\tableau*[scY]{5\cr 4&6\cr2&3&4&4&6\cr 1&2&2&2&3&4&4&6 }}} \quad
{\tiny{\tableau*[scY]{6\cr 4&5\cr 2&3&4&4&5\cr 1&2&2&2&3&4&4&5 }}} \quad
{\tiny{\tableau*[scY]{4\cr 3&6\cr 2&4&4&5&6\cr 1&2&2&2&4&4&5&6 }}}
\end{equation}
\end{example}

\medskip

\subsection{Standardizing and deleting a letter from $k$-tableaux}
As with the standard $k$-tableaux, we shall prove that deleting 
some letter from a $k$-tableau produces another $k$-tableau.
To this end, we present a method for constructing a standard
$k$-tableau from a given $k$-tableau of the same shape.

\begin{definition}  
For $\alpha$ a composition of $m$ and $T \in \mathcal T_{\alpha}^{k}(\lambda)$, 
define $St(T)$ by the iterative process:
\begin{quote}
\noindent If $a$ is the biggest letter of $T$, let
$i$ denote the $k+1$-residue of the rightmost cell in $T$ 
that contains $a$.  Construct a tableau $\bar T$ by replacing 
each occurrence of letter $a$ with $k+1$-residue $i$ by the 
letter $m$.  Now let $a$ denote the biggest letter (smaller 
than $m$) in $\bar T$ and $i$ the $k+1$-residue of the rightmost 
cell in $\bar T$ that contains $a$. Again construct a new tableau by 
replacing each occurrence of letter $a$ with $k+1$-residue $i$ 
by the letter $m-1$.  $St(T)$ is the tableau obtained by iterating
this process until each collection of repeated letters forms only one $k+1$-string.
That is, $St(T)\in\mathcal T_{1^m}^k(\lambda)$.
\end{quote}
\end{definition}

\begin{example}
Given a $k$-tableau $T\in\mathcal T_{(1,3,1,2,1,1)}^{3}
(3,3,2,1)$ of shape
$\core\left(3,3,2,1\right)=(8,5,2,1)$:
\begin{equation}
\nonumber
T\;=\; {\tiny{\tableau*[scY]{5\cr 4&6\cr2&3&4&4&6\cr 1&2&2&2&3&4&4&6 }}} 
\quad
\;\; \text{\rm 4-residues} \;=\;
{\tiny{\tableau*[scY]{1\cr 2&3\cr3&0&1&2&3\cr 0&1&2&3&0&1&2&3 }}} 
\qquad \qquad \qquad \qquad
\qquad \qquad \qquad \qquad
\end{equation}
\begin{equation}
\nonumber
\text{\rm Every letter $a=6$ of residue $i=3$ is replaced by $m=9$:}
\;\;
{\tiny{\tableau*[scY]{5\cr 4&9\cr2&3&4&4&9\cr 1&2&2&2&3&4&4&9 }}} 
\qquad \qquad \qquad
\qquad \qquad \qquad
\end{equation}
\begin{equation}
\nonumber
\text{\rm Then letters $a=5$ of residue $i=1$ are replaced by $m=8$:}
\;\;
{\tiny{\tableau*[scY]{8\cr 4&9\cr2&3&4&4&9\cr 1&2&2&2&3&4&4&9 }}} 
\qquad \qquad \qquad
\qquad \qquad \qquad
\end{equation}
\begin{equation}
\nonumber
\text{\rm Then letters $a=4$ of residue $i=2$ are replaced by $m=7$:}
\;\;
{\tiny{\tableau*[scY]{8\cr 7&9\cr2&3&4&7&9\cr 1&2&2&2&3&4&7&9 }}} 
\qquad \qquad \qquad
\qquad \qquad \qquad
\end{equation}
\begin{equation}
\nonumber
\text{\rm Then letters $a=4$ of residue $i=1$ are replaced by $m=6$:}
\;\;
{\tiny{\tableau*[scY]{8\cr 7&9\cr2&3&6&7&9\cr 1&2&2&2&3&6&7&9 }}} 
\qquad \qquad \qquad
\qquad \qquad \qquad
\end{equation}

\begin{equation}
\text{\rm Similarly,}\;
{{a=3, i=0\atop\longrightarrow}\atop m=5}
\; \;
\; \;
{\tiny{\tableau*[scY]{8\cr 7&9\cr2&5&6&7&9\cr 1&2&2&2&5&6&7&9 }}} 
\; \;
\; \;
{{a=2, i=3\atop\longrightarrow}\atop m=4}
\; \;
\; \;
{\tiny{\tableau*[scY]{8\cr 7&9\cr4&5&6&7&9\cr 1&2&2&4&5&6&7&9 }}} 
\; \;
\; \;
{{a=2, i=2\atop\longrightarrow}\atop m=3}
\; \;
\; \;
{\tiny{\tableau*[scY]{8\cr 7&9\cr4&5&6&7&9\cr 1&2&3&4&5&6&7&9 }}} 
\qquad \qquad \qquad
\qquad \qquad \qquad
\nonumber
\end{equation}

\noindent
Once the tableau is standard, the step $a=2,i=1,m=2$ followed 
by $a=1,i=0,m=1$ does not change the tableau. 
\end{example}

\begin{proposition} 
\label{prop70} 
Let $T \in \mathcal T_{(\alpha_1,\ldots,\alpha_m)}^k(\lambda)$.  
The tableau obtained by deleting the letter $m$ from $T$
belongs $\mathcal T_{(\alpha_1,\dots,\alpha_{m-1})}^k(\mu)$ 
for some $\mu \preceq \lambda$ with $|\lambda|-|\mu|=\alpha_m$.
\end{proposition}
\begin{proof}
Let $\hat T$ denote the tableau obtained by deleting the letter $m$ 
from $T$.  Since Conditions~(i) and (ii) of a $k$-tableau
clearly hold for $\hat T$, it suffices to show that the shape of $\hat T$ is
given by ${\core}(\mu)$ for some $\mu\preceq\lambda$.
To this end, consider $St(T)$, the standard $k$-tableau of 
shape ${\core}(\lambda)$ associated to $T$.
Deleting the largest letter from $St(T)$ gives a $k$-standard tableau 
of shape $s_{i_{\aa_m}}\bigl({\core}(\lambda) \bigr)$ 
by Proposition~\ref{si=hatsi}.
By iteration, removing the largest $\aa_m$ letters from $St(T)$ 
gives a standard $k$-tableau $\bar T$ of shape
$s_{i_{1}}\cdots s_{i_{\alpha_m}}\bigl({\core}(\lambda)\bigr)$,
where ${i_{1}},\dots,{i_{\alpha_m}}$ are respectively the $k+1$-residues 
of the $\aa_m$ largest letters in $St(T)$.  Since $\hat T$ has the 
same shape as $\bar T$, 
$\hat T\in\mathcal T_{(\alpha_1,\dots,\alpha_{m-1})}^k(\mu)$ where 
${\core}(\mu)=s_{i_{1}}\cdots s_{i_{\alpha_m}}\bigl({\core}
(\lambda)\bigr)$.  Further, 
$\mu\preceq\lambda$ by Corollary~\ref{corokhook}
and $|\lambda|-|\mu|=\alpha_m$ by Proposition~\ref{remarkone}.
\end{proof}

It is known (eg. \cite{[Fu]}) that there are no semi-standard tableaux 
of shape $\lambda$ and evaluation $\mu$ when $\lambda \ntrianglerighteq \mu$
in dominance order.  We have found that this is also true for the 
$k$-tableaux.

\begin{remark}
There are no $k$-tableaux in $\mathcal T_\mu^{k}(\lambda)$ when 
$\ell(\mu)<\ell(\lambda)$ since any element of $\mathcal T_\mu^{k}(\lambda)$ 
has height $\ell(\lambda)$, has only $\ell(\mu)$ distinct letters, and
must be strictly increasing in columns.
\label{shapecontent}
\end{remark}

\begin{theorem}
\label{triangu}
There are no semi-standard $k$-tableaux in $\mathcal T_\mu^{k}(\lambda)$
when $\lambda \ntrianglerighteq \mu$.  Further, there is exactly 
one when $\lambda=\mu$.
\end{theorem}
\begin{proof}
Consider $\lambda,\mu\in\mathcal P_k$ with $|\lambda|=|\mu|$.
We shall proceed by induction on the length of $\mu$.  
A $k$-tableau of evaluation $\mu=(\mu_1)$ must be of shape
${\core}(\lambda)$ where $\ell(\lambda)\leq \ell(\mu)$
by Remark~\ref{shapecontent}.  Therefore, $\lambda=(\mu_1)$ and the 
claim holds.  Assume the assertion holds when $\ell(\mu)<N$.  

Consider 
$\mu=(\mu_1,\ldots,\mu_N)$ and $\lambda=(\lambda_1,\ldots,\lambda_{\bar N})$
with $\lambda \ntrianglerighteq \mu$. That is,
$\mu_1+\cdots+\mu_j>\lambda_1+\cdots+\lambda_j$ for some 
$j\leq N$.  Suppose by contradiction that there is some 
$T\in\mathcal T_\mu^{k}(\lambda)$.  The previous proposition
implies that removing the letter $N$ from $T$ results in a $k$-tableau 
$\bar T\in\mathcal T_{\bar\mu=(\mu_1,\ldots,\mu_{N-1})}^k(\bar\lambda)$ 
where $\bar\lambda\preceq\lambda$.  Thus,
$\core(\bar\lambda)= s_{i_1}\cdots  
s_{i_{\mu_N}}\bigl(\core(\lambda)\bigr)$
for some $i_1,\ldots,i_{\mu_N}$ by Corollary~\ref{corokhook}.
Since $\ell(\bar\mu)<N$, the induction hypothesis implies 
$\bar\mu\unlhd\bar\lambda$.  
Therefore $\mu_1+\cdots+\mu_r\leq\bar\lambda_1+\ldots+\bar\lambda_r$ 
for all $r\leq N-1$.  Further, $\bar\lambda_i\leq\lambda_i$ 
by Proposition~\ref{remarkone} since the $s_{i_{j}}$ act by deleting 
removable corners starting with $\core(\lambda)$, and thus
$\mu_1+\cdots+\mu_r\leq\lambda_1+\ldots+\lambda_r$ for all $r\leq N-1$.
Therefore, $\mu_1+\cdots+\mu_N>\lambda_1+\cdots+\lambda_N$ given
$\lambda\ntrianglerighteq\mu$.
However, since $|\lambda|=|\mu|$, $|\mu|>\lambda_1+\cdots+\lambda_N$
implies $\ell(\lambda)>\ell(\mu)$.  
We thus reach a contradiction by
Remark~\ref{shapecontent}.

To see that there is exactly one $k$-tableau 
$T\in\mathcal T_{\lambda}^k(\lambda)$, we shall first show 
by induction that there can be no more than one such 
tableau for $\lambda=(\lambda_1,\ldots,\lambda_N)$.
Delete the letter $N$ from $T$ to obtain a $k$-tableau 
$\bar T\in\mathcal T_{\hat\lambda=(\lambda_1,\ldots,\lambda_{N-1})}^k
(\bar\lambda)$ where
$\core(\bar\lambda)=s_{i_1}\cdots s_{i_{\lambda_N}}\bigl(\core(\lambda)\bigr)$.
Remark~\ref{shapecontent} implies that
$\ell(\bar\lambda)\leq \ell(\hat\lambda)=N-1$.
Since exactly $\lambda_N$ cells were removed from 
$\lambda$ to obtain $\bar\lambda$, and the length 
of $\lambda$ was decreased by at least one, the
${i_1},\ldots,{i_{\lambda_N}}$ are uniquely determined 
and correspond to the $k+1$-residues in the top row of 
${\core}(\lambda)$.  Thus, for two distinct $k$-tableaux in 
$\mathcal T_{\lambda}^k(\lambda)$ to exist, two distinct
$k$-tableaux in $\mathcal T^k_{\hat\lambda}(\hat \lambda)$ 
are necessary.  By induction this is a contradiction.  

We prove that there is in fact always a 
$k$-tableau $T\in\mathcal T_{\lambda}^k(\lambda)$ by
construction: start with
unique $k$-tableau of shape and evaluation $(\lambda_1)$.  For $j\geq 1$, 
let $\lambda^{(j)}=(\lambda_1,\ldots,\lambda_j)$ and consider the 
$k$-tableau of shape ${\core}(\lambda^{(j)})$ and evaluation 
$\lambda^{(j)}$.  Add the letter $j+1$ in all positions 
$s_{i_{\lambda_{j+1}}}\cdots s_{i_1}\bigl({\core}(\lambda^{(j)})\bigr)
\big/{\core}(\lambda^{(j)})$ where $i_\ell$ is the $k+1$-residue 
of the square $(j+1,\ell)$ of $\lambda^{(j)}$ for 
$\ell=1,\ldots,\lambda_{j+1}$.
\end{proof}

\subsection{Bijection: generalized $k$-tableaux and chains in the $k$-lattice}
A rule for expanding the product of a $k$-Schur function with the 
homogeneous function $h_\ell$ (for $\ell\leq k$) in terms of $k$-Schur 
functions was conjectured in \cite{[LM1]}.  We introduce certain 
sequences of partitions based on this generalized Pieri rule and
show their connection to the semi-standard $k$-tableaux.
The connection with symmetric functions is then discussed in
\S\ref{symfu}.

Recall from the introduction that a pair of $k$-bounded partitions 
$\la,\mu$ is ``{\it $r$-admissible}'' if and only if $\la/\mu$ and  
$\la^{\om_k}/\mu^{\om_k}$ are respectively horizontal and vertical 
$r$-strips.  For composition $\aa$, a sequence of partitions
$\left(\la^{(0)},\la^{(1)},\cdots,\la^{(r)}\right)$
is ``$\alpha$-admissible" if $ \la^{(j)},\la^{(j-1)}$ is 
a $\alpha_j$-admissible pair for all $j$.  This given, since 
Corollary~\ref{coronouveau} implies that if 
$\lambda^{(j)},\lambda^{(j-1)}$ is $\alpha_j$-admissible
then $\lambda^{(j-1)}\preceq\lambda^{(j)}$,
we have that any $\alpha$-admissible
sequence must be a chain in the $k$-Young
lattice.  We are interested in the set of chains:

\begin{definition}  
\label{coveringdef} 
For any composition $\aa$, let
$$
\mathcal D^k_{\aa}(\lambda)=\left\{ 
(\emptyset=\lambda^{(0)},\ldots,\lambda^{(r)}=\lambda) 
\;\;\text{that are}\;\; \alpha\text{-admissible}
\right\}\,.
$$
\end{definition}

We now give a bijection between the set of chains
in $\mathcal D^k_{\aa}(\lambda)$ and the tableaux in 
$\mathcal T^k_{\aa}(\lambda)$.

\begin{definition} 
For any $P=(\lambda^{(0)},\lambda^{(1)},\dots,\lambda^{(m)})
\in {\mathcal D}_{\aa}^k(\lambda)$, let $\Gamma(P )$
be the tableau of shape ${\core}(\lambda)$ 
where letter $j$ fills cells in positions 
${\core}(\lambda^{(j)})/{\core}(\lambda^{(j-1)})$, for $j=1,\dots,m$.
\end{definition}

\medskip

\begin{proposition}  \label{propermapout}
If  $ P \in {\mathcal D}_{\aa}^k(\lambda)$, then   
$\Gamma(P) \in {\mathcal T}_{\aa}^k(\lambda)$. 
\end{proposition}
\begin{proof}
If $P=(\lambda^{(0)},\lambda^{(1)},\dots,\lambda^{(m)}) 
\in \mathcal D_{\aa}^k(\lambda)$ then
$\Gamma(P)$ has the shape of the $k+1$-core ${\core}(\lambda)$.  It thus
suffices to prove that
$\Gamma(P)$ is column strict and has $\alpha_j$  
distinct $k+1$-residues that are filled with the letter $j$.
Since $\lambda^{(j)},\lambda^{(j-1)}$ is $\alpha_j$-admissible
by definition of $\mathcal D_{\aa}^k(\lambda)$,
Theorem~\ref{coro67} implies that
$s_{i_1}\cdots s_{i_{\aa_j}}\bigl({\core}(\lambda^{(j-1)}) \bigr)
={\core}(\lambda^{(j)})$ for some collection of distinct integers
$i_1,\ldots,i_{\aa_j}$ and Proposition~\ref{propo64} implies that
${\core}(\lambda^{(j)})/{\core}(\lambda^{(j-1)})$ is a horizontal strip. 
$\Gamma(P)$ is thus column strict since the letter $j$ lies only 
in a horizontal strip.  Further, given that each of the $\aa_j$ operators $s_{i_t}$ 
adds addable corners of residue $i_t$, the letter $j$ occupies 
$\aa_j$ distinct $k+1$-residues since 
$i_1,\ldots,i_{\aa_j}$ are distinct.
\end{proof}

\medskip

\begin{definition}
For a $k$-tableau $T\in\mathcal T_\aa^{k}(\lambda)$
with $\alpha=(\alpha_1,\dots,\alpha_m)$,
let $\bar\Gamma(T)=(\lambda^{(0)},\ldots,\lambda^{(m)})$,
where $\core(\lambda^{(i)})$ is the shape of the tableau
obtained by deleting the letters $i+1,\dots,m$ from $T$.
\end{definition}

\medskip

\begin{proposition}  \label{propomapin}
If $T \in \mathcal T_{\alpha}^k(\lambda)$, 
then $\bar \Gamma (T)\in\mathcal D_\aa^k(\lambda)$. 
\end{proposition}
\begin{proof}  
Letting $\bar \Gamma (T)=(\lambda^{(0)},\dots,\lambda^{(m)} )$,
the definition of $\bar\Gamma$ implies that the tableaux $T^i$ 
obtained by deleting letters $i+1,\dots,m$ from $T$ 
have corresponding shapes $\core(\lambda^{(i)})$.  
By Proposition~\ref{prop70}, the $T^i$ are $k$-tableaux.  
In particular, $T^i$ has strictly increasing columns.  
Thus since $T^{i-1}$ is obtained by deleting letter $i$ from $T^i$,
${\core}(\lambda^{(i)})/{\core}(\lambda^{(i-1)})$ is a 
horizontal strip and further, by Proposition~\ref{prop70},
$\lambda^{(i-1)} \preceq \lambda^{(i)}$  with 
$|\lambda^{(i)}|-|\lambda^{(i-1)}|=\alpha_i$.  
Property~\ref{contains} then implies that 
$\lambda^{(i-1)} \subseteq \lambda^{(i)}$ and 
$(\lambda^{(i-1)})^{\omega_k}\subseteq(\lambda^{(i)})^{\omega_k}$.  
Therefore $\lambda^{(i)},\lambda^{(i-1)}$ are $\alpha_i$-admissible
by Theorem~\ref{lemhoriverti} and we have
that $\bar \Gamma (T) \in \mathcal D_\aa^k(\lambda)$.
\end{proof}

\begin{theorem} \label{bijecgenn}
$\Gamma$  is a bijection between
$\mathcal T_\aa^k(\lambda)$ and $\mathcal D_\aa^k(\lambda)$,
with $\Gamma^{-1}=\bar \Gamma$.
\end{theorem}
\begin{proof}
Given Propositions~\ref{propermapout} and \ref{propomapin},
we only have to show that if $P \in \mathcal D_\alpha^k(\lambda)$ and
$T \in \mathcal T_\alpha^k(\lambda)$, 
then $\Gamma (\bar \Gamma(T))=T$ and
$\bar \Gamma (\Gamma(P))=P$.  This follows from the 
same deleting-filling letter argument given in the proof of  
Theorem~\ref{firstbijection}.
\end{proof}

\medskip

\section{Symmetric functions and $k$-tableaux}
\label{symfu}

Refer to \cite{[M]} for details on symmetric functions and
Macdonald polynomials. 
Here, we are interested in the study of the $q,t$-Kostka 
polynomials $K_{\mu \lambda}(q,t) \in \mathbb N[q,t]$.
These polynomials arise as expansion coefficients for the 
Macdonald polynomials $J_\lambda[X;q,t]$ in terms of a basis dual to 
the monomial basis with respect to the Hall-Littlewood scalar product.  
As in the introduction, we use the modification of $J_\lambda[X;q,t]$ whose expansion 
coefficients in terms of Schur functions are the $q,t$-Kostka 
coefficients: 
\begin{equation}
\label{schurexp}
H_{\lambda}[X;q,t]=\sum_{\mu} K_{\mu \lambda}(q,t)\, 
s_{\mu}[X] \, .
\end{equation}
The $q,t$-Kostka coefficients also have a representation theoretic 
interpretation \cite{[GH]}, from which they were shown \cite{[H]} to 
lie in $\mathbb N[q,t]$.
Since $J_{\lambda}[X;q,t]$ reduces to the Hall-Littlewood
polynomial $Q_\lambda[X;t]$ when $q=0$, we obtain a modification of the 
Hall-Littlewood polynomials by taking:
\begin{equation}
H_{\lambda}[X;t]=H_{\lambda}[X;0,t]=
\sum_{\mu\unrhd \lambda} K_{\mu \lambda}(t)\, s_{\mu}[X] \, ,
\label{HallinS}
\end{equation}
with the coefficients $K_{\mu \lambda}(t)\in\mathbb N[t]$ known 
as Kostka-Foulkes polynomials.  We can then obtain the homogeneous 
symmetric functions by letting $t=1$:
\begin{equation} \label{kostka}
h_{\lambda}[X]=H_{\lambda}[X;1]= 
\sum_{\mu\unrhd\lambda} K_{\mu \lambda}\, s_{\mu}[X] \, ,
\end{equation}
where $K_{\mu \lambda} \in \mathbb N$ are the Kostka numbers.  

\medskip

Recent work in the theory of symmetric functions has led to a new
approach in the study of the $q,t$-Kostka polynomials.
The underlying mechanism for this approach relies on a family of 
polynomials that appear to have a remarkable kinship with the 
classical Schur functions \cite{[LLM],[LM1],[LM2]}.
More precisely,
consider the filtration $\Lambda_t^{(1)}\subseteq \Lambda_t^{(2)}\subseteq
\cdots \subseteq\Lambda_t^{(\infty)}=\Lambda$,
given by linear spans of Hall-Littlewood polynomials indexed by $k$-bounded 
partitions.  That is,
\def \LAt {\Lambda_t^{(k)}}
$$
\Lambda^{(k)}_t = 
{\cal L} \{H_{\lambda}[X;t] \}_{\lambda;\lambda_1 \leq k}\, , 
\qquad k=1,2,3,\dots \, .
$$ 
A family of symmetric functions called the $k$-Schur functions, 
$s_\la^{(k)}[X;t]$, was introduced in \cite{[LM1]} (these functions 
are conjectured to be precisely the polynomials defined using 
tableaux in \cite{[LLM]}).  It was shown that the $k$-Schur functions
form a basis for $\LAt$ and that, for $\lambda$ a $k$-bounded partition, 
\begin{equation}
H_{\lambda}[X;q,t \, ] \;= \;\sum_{\mu;\mu_1\leq k} 
K_{\mu \lambda}^{(k)}(q,t) \, s_{\mu}^{(k)}[X;t\, ] \, , \qquad
K_{\mu \lambda}^{(k)}(q,t) \in \Z[q,t] \, ,
\end{equation}
and 
\begin{equation}
\label{hallkschur}
H_{\lambda}[X;t \, ]\; = \;
s_\lambda^{(k)}[X;t\,] \;+\;
\sum_{\mu;\mu_1\leq k\atop \mu>_D\lambda} 
K_{\mu \lambda}^{(k)}(0,t) \, s_{\mu}^{(k)}[X;t\, ] \, , \qquad
K_{\mu \lambda}^{(k)}(0,t) \in \Z[t] \, .
\end{equation}
The study of the $k$-Schur functions is motivated in part by
the conjecture \cite{[LLM],[LM1]} that the expansion
coefficients actually lie in $\mathbb N[q,t]$.  That is,
\begin{equation}
K_{\mu \lambda}^{(k)}(q,t) \in \N[q,t] \, .
\label{macconjec} 
\end{equation}
Since it was shown that $s_\lambda^{(k)}[X;t]=s_\lambda[X]$
for $k$ larger than the hook-length of $\lambda$,
this conjecture generalizes Eq.~\eqref{schurexp}.
Also, there is evidence to support that
$K_{\mu\lambda}(q,t)-K_{\mu\lambda}^{(k)}(q,t)\in\mathbb N[q,t]$,
suggesting that the $k$-Schur expansion coefficients 
are simpler than the $q,t$-Kostka polynomials.

\medskip

The preceding developments on the $k$-lattice can be applied 
to the study of the generalized $q,t$-Kostka coefficients
as follows: the $k$-Schur functions appear to obey a generalization of the 
Pieri rule on Schur functions.  To be precise, it was conjectured 
in \cite{[LLM],[LM1]} that for the complete symmetric function $h_{\ell}[X]$,
\begin{equation} 
\label{pierih}
h_{\ell}[X]\,s_{\lambda}^{(k)}[X;1] =
\sum_{\mu \in E_{\lambda,\ell}^{(k)}}s_{\mu}^{(k)}[X;1] \, ,
\end{equation}
where 
\begin{equation}E_{\lambda,\ell}^{(k)}=\bigl \{ \mu \, | \, \mu/\lambda 
{\rm{~is~a~horizontal~}}\ell{\rm{-strip~and ~}} 
\mu^{\omega_k}/\lambda^{\omega_k} 
{\rm{~is~a~vertical~}}\ell{\rm{-strip}} \bigr \}\,.
\label{piericover}
\end{equation}
Iteration, from $s_{\emptyset}^{(k)}[X;1]=1$, then yields that the 
expansion of $h_{\lambda_1}[X] h_{\lambda_2}[X] \cdots$ satisfies
\begin{equation} \label{hkschur}
h_{\lambda}[X] = \sum_{\mu} K_{\mu \lambda}^{(k)} \, s_{\mu}^{(k)}[X;1] \, ,
\end{equation}
where $K_{\mu \lambda}^{(k)}$ is a nonnegative integer reducing to
the usual Kostka number $K_{\mu \lambda}$ when $k$ is large since
$s_\lambda^{(k)}[X;t]=s_\lambda[X]$ in this case.
The definition of $E_{\lambda,\ell}^{(k)}$ in the $k$-Pieri expansion 
thus reveals the motivation behind the set of chains given 
in Definition~\ref{coveringdef}. This connection implies
that
\begin{equation}
\nonumber
K_{\mu \lambda}^{(k)} 
\, = \,\text{the number of chains 
of the $k$-lattice in $\mathcal D_{\lambda}^k(\mu)$}
\,.
\end{equation}
Equivalently, using the bijection between 
chains in $\mathcal D_\la^k(\mu)$ and $\mathcal T_\la^k(\mu)$
given in Theorem~\ref{bijecgenn}, we have 
\begin{equation}
\nonumber
K_{\mu \lambda}^{(k)}
\,=\,
\text{
the number of $k$-tableaux of shape
${\core}(\mu)$ and evaluation $\la$}
\,.
\end{equation}
Although this combinatorial interpretation relies on the conjectured 
Pieri rule \eqref{pierih}, it was proven in \cite{[LM1]} that
the $k$-Schur functions are unitriangularly related to the
homogeneous symmetric functions.  That is,
$K_{\lambda\mu}^{(k)}=0$ when $\mu \ntrianglerighteq \lambda$ and 
$K_{\lambda\lambda}^{(k)}=1$.  Therefore,
Theorem~\ref{triangu} implies that the number of 
$k$-tableaux does correspond to $K_{\lambda\mu}^{(k)}$
in these cases.

\medskip

More generally, note that letting $q=0$ in Eq.~\eqref{macconjec} 
gives that the coefficients in Hall-Littlewood expansion 
Eq.~\eqref{hallkschur} satisfy
$K^{(k)}_{\mu\lambda}(0,t)\in\mathbb N[t]$.
However, since $H_{\lambda}[X;1]=h_{\lambda}[X]$,
we have that $K_{\mu \lambda}^{(k)}(0,1)=K_{\mu \lambda}^{(k)}$
from Eq.~\eqref{hkschur}.  Therefore, since it appears that
$K_{\mu \lambda}^{(k)}$ counts the number of semi-standard $k$-tableaux
in $\mathcal T_{\lambda}^k(\mu)$,
it is suggested that there exists a $t$-statistic on such 
$k$-tableaux giving a combinatorial 
interpretation for 
the generalized Kostka-Foulkes $K_{\mu \lambda}^{(k)}(0,t)$.

\medskip

Alternatively, $H_{\lambda}[X;1,1]=h_{1^n}[X]$ for $\lambda\vdash n$
implies that $K_{\mu \lambda}^{(k)}(1,1)=K_{\mu \, 1^n}^{(k)}$
by Eq.~\eqref{hkschur}.  This lends support to the idea that
a $q,t$-statistic on the standard $k$-tableaux that would account 
for the apparently positive coefficients $K_{\mu \lambda}^{(k)}(q,t)$ in 
Eq.~\eqref{macconjec}.  That is,
\begin{equation}
\nonumber
K_{\mu \lambda}^{(k)}(1,1)\,=\, \text{
the number of standard $k$-tableaux of shape 
${\core}(\mu)$}
\,.
\end{equation}

\medskip

Equivalently, our bijection between affine permutations and standard $k$-tableaux
suggests there may be a $q,t$-statistic on reduced words
that would account for the positivity:
\begin{equation}
\nonumber
K_{\mu \lambda}^{(k)}(1,1)\,=\, \text{
the number of reduced words of $\sigma_\mu\in 
\tilde S_{k+1}/S_{k+1}$ }
\,.
\end{equation}

We mention one final consequence of the $k$-Pieri rule.
For $\lambda$ a partition of length $n$,
the product $h_{\lambda_1}\cdots h_{\lambda_n}$ giving
$h_\lambda$ can be written in any order since the functions commute.
Therefore, 
\begin{equation} 
h_{\alpha}[X] = \sum_{\mu} K_{\mu \lambda}^{(k)} \, s_{\mu}^{(k)}[X;1] \, .
\end{equation}
for any reordering $\alpha$ of the entries of $\lambda$.
Therefore, $K_{\mu \lambda}^{(k)}$ is also 
the number of chains in $\mathcal D_{\alpha}^k(\mu)$.  Equivalently, 
$K_{\mu \lambda}^{(k)}$ is the number of $k$-tableaux in 
$\mathcal T_{\alpha}^k(\mu)$.  Thus, conjecture \eqref{pierih} implies:

\begin{quote}
If $\alpha$ is a rearrangement of $\lambda$,
then $\big|\mathcal T_{\alpha}^k(\mu)\big|=
\big|\mathcal T_{\lambda}^k(\mu)\big|$.  Equivalently,
the number of $k$-tableaux in $\mathcal T_{\alpha}^k(\mu)$
equals the number of $k$-tableaux in
$\mathcal T_{\lambda}^k(\mu)$.
\end{quote}

This conjecture suggests that there is a generalization of the
Bender-Knuth involution on semi-standard tableaux that permutes 
the evaluation of $k$-tableaux accounting for this phenomenon.
See \cite{[LMW]} for this new involution and thus the proof of this
conjecture.

\noindent
{\bf Acknowledgments.}
{\it 
We would like to thank Michelle Wachs for helpful
and revealing correspondence.  
We are thankful 
to Ira Gessel for his tableaux macros.}

\end{document}